\documentclass{amsart}[10pt]
\usepackage{enumerate, amsfonts,pstcol,epsfig,amssymb, amsmath, amsthm}

\newlength{\tabwidth}
\newlength{\tabheight}
\newlength{\tabrule}
\newlength{\tabwidthx}
\newlength{\tabheightx}

\def\gentabbox#1#2#3#4{\vbox to \tabheight{\setlength{\tabrule}{#3}%
  \setlength{\tabwidthx}{#1\tabwidth}\addtolength{\tabwidthx}{\tabrule}%

\setlength{\tabheightx}{#2\tabheight}\addtolength{\tabheightx}{-\tabheight}%
  \hbox to #1\tabwidth{%
 \hspace{-0.5\tabrule}\rule{\tabrule}{#2\tabheight}\hspace{-\tabrule}%
    \vbox to #2\tabheight{\hsize=\tabwidthx%
      \vspace{-0.5\tabrule}\hrule width\tabwidthx height\tabrule%
      \vspace{-0.5\tabrule}\vfil%
      \hbox to \tabwidthx{\hss#4\hss}%
        \vfil\vspace{-0.5\tabrule}%
      \hrule width\tabwidthx height\tabrule\vspace{-0.5\tabrule}}%
 \hspace{-\tabrule}\rule{\tabrule}{#2\tabheight}\hspace{-0.5\tabrule}}%
  \vspace{-\tabheightx}}}
\def\genblankbox#1#2{\vbox to \tabheight{\vfil\hbox to
#1\tabwidth{\hfil}}}
\def\tabbox#1#2#3{\gentabbox{#1}{#2}{0.4pt}{\strut #3}}

\catcode`\:=13 \catcode`\.=13 \catcode`\;=13
\catcode`\>=13 \catcode`\^=13
\def:#1\\{\hbox{$#1$}}
\def.#1{\tabbox{1}{1}{$#1$}}
\def>#1{\tabbox{2}{1}{$#1$}}
\def^#1{\tabbox{1}{2}{$#1$}}
\def;{\genblankbox{1}{1}\relax}
\catcode`\:=12 \catcode`\.=12 \catcode`\;=12
\catcode`\>=12 \catcode`\^=7

\newenvironment{tableau}{\bgroup\catcode`\:=13 \catcode`\.=13
  \catcode`\;=13 \catcode`\>=13 \catcode`\^=13
  \setlength{\tabheight}{3ex}\setlength{\tabwidth}{3ex}%
  \def\b##1##2##3{\gentabbox{##1}{##2}{1.2pt}{\vbox{##3}}}%
  \def\n##1##2##3{\gentabbox{##1}{##2}{0.4pt}{\vbox{##3}}}%
  \vbox\bgroup\offinterlineskip}{\egroup\egroup}

\newtheorem{theorem}{Theorem}[section]
\newtheorem{corollary}[theorem]{Corollary}
\newtheorem{lemma}[theorem]{Lemma}
\newtheorem{proposition}[theorem]{Proposition}

\newtheorem{fact}{Fact}[section]
\newtheorem*{theorem:main}{Theorem \ref{theorem:generators}}
\newtheorem*{lemma:technical2}{Lemma \ref{lemma:technical2}}
\newtheorem*{prop:amongspecial}{Proposition \ref{prop:amongspecial}}
\newtheorem*{prop:makespecial}{Proposition \ref{prop:makespecial}}

\newtheorem*{conjecture*}{Conjecture}

\theoremstyle{definition}
\newtheorem{definition}[theorem]{Definition}

\theoremstyle{remark}

\newtheorem{example}[theorem]{Example}
\newtheorem{remark}[fact]{Remark}

\begin{document}
\title{Knuth Relations for the Hyperoctahedral Groups}
\author{Thomas Pietraho}
\email{tpietrah@bowdoin.edu} 
\address{Department of Mathematics\\Bowdoin College\\Brunswick,
Maine 04011}
 \maketitle

\begin{abstract}
C.~Bonnaf{\'e}, M.~Geck, L.~Iancu, and T.~Lam
have conjectured a description of one-sided cells in unequal parameter Hecke algebras of type $B$ which is based on domino tableaux of arbitrary rank.  In the integer case, this generalizes the work of D.~Garfinkle whose methods we adapt to construct a family of operators which generate the conjectured combinatorial description.
\end{abstract}
\section{Introduction}

In \cite{garfinkle3}, D.~Garfinkle classified the primitive spectrum of the universal enveloping algebra for a complex semisimple Lie algebra in types $B$ and $C$.  By using annihilators of highest weight modules, this problem is reduced to studying equivalence classes in the corresponding Weyl group $W_n$.  The existence of a Robinson-Schensted bijection between elements of $W_n$ and same shape pairs of standard domino tableaux with $n$ dominos \cite{garfinkle1} turns this into an essentially combinatorial problem.  In fact,  Garfinkle's classification shows that two elements in $W_n$ are equivalent iff their left domino tableaux are related by moving through an open cycle, a certain combinatorial operation. The key step of this classification was achieved by studying the action of the wall-crossing operators arising from the general $\tau$-invariant, as defined in  \cite{vogan:tau}, which are shown to be generators for both equivalences.

When interpreted in terms of cells of the equal parameter Hecke algebra of type $B$, the above work takes on a new meaning.
In \cite{lusztig:unequal}, G.~Lusztig extended the notion of a Hecke algebra associated to a Coxeter group by introducing a weight function to its defining relations.  For Weyl groups of type $B$, cell theory then depends on an additional parameter $s$, which reduces to the equal parameter case when $s=1$;  Garfinkle's work classified cells in  exactly this latter setting.   As observed in \cite{vanleeuwen:rank}, Garfinkle's bijection also admits an extension
$$G_r: W_n \rightarrow SDT_r(n) \times SDT_r(n)$$
to a bijection between $W_n$ and same shape pairs of standard domino tableaux of rank $r$.  It is reasonable to hope that the above two parameters can be linked and that
a similar classification of cells is possible in this more general case. In fact:

\begin{conjecture*}[\cite{bgil}] When $s$ is a positive integer, two elements of $W_n$ lie in the same left cell if and only if their left domino tableaux of rank $s-1$ are related
by moving through an open cycle.
\end{conjecture*}

This is of course true when $s=1$ and has been verified in the asymptotic case $s\geq n$, when the bijection $G_r$ degenerates to the generalized Robinson-Schensted correspondence of \cite{stanton:white}, see \cite{bonnafe:iancu}.  Among finite Coxeter groups, cells in unequal parameter Hecke algebras have been classified in the dihedral groups and type $F4$ by  Geck and Pfeiffer \cite{geck:pfeiffer}, Geck \cite{geck:leftcells}, and Lusztig \cite{lusztig:unequal}.  Only the problem of their classification in type $B$ remains.

We will say that two elements of the hyperoctahedral group $W_n$ are in the same irreducible combinatorial left cell of rank $r$ if they share the same left domino tableau under the Robinson-Schensted map $G_r$, and in the same reducible combinatorial left cell of rank $r$ if their rank $r$ left domino tableaux are related by moving through an open cycle. The previous conjecture can be restated as:

\begin{conjecture*}[\cite{bgil}] When $s$ is a positive integer, left cells in $W_n$ coincide with reducible combinatorial left cells of rank $s-1$.
\end{conjecture*}

 Inspired by Garfinkle's approach in the equal parameter case, the main goal of this paper is to construct a set of generators for the reducible combinatorial left cells in arbitrary rank which draws on the notion of the generalized $\tau$-invariant used in the equal parameter case.
Such a set $\Lambda^{r+1}$ is constructed in Section \ref{section:knuth}.  The main theorem then can be stated as:

\begin{theorem:main} The family of operators $\Lambda^{r+1}$ generates the reducible combinatorial left cells of rank $r$.  More precisely, given $w$ and $v \in W_n$ whose rank $r$ left domino tableaux  differ only by moving through a set of non-core open cycles, there is a sequence of operators in $\Lambda^{r+1}$ sending $w$ to $v$.
\end{theorem:main}

This result falls into a family of similar theorems on generating sets for equivalence classes of standard tableaux, Garfinkle's original not withstanding.
In type $A$, a family of operators generating a similar equivalence for the symmetric group on standard Young tableaux appears in \cite{knuth:art3} and is known as the set of Knuth relations.  Domino tableaux whose shapes are also partitions of nilpotent orbits in types $B$, $C$, or $D$ correspond to the so-called orbital varieties.  Their classification has been carried out by W.~M.~McGovern in \cite{mcgovern:ssmap} by relying on a similar set of generators found in \cite{hopkins:domino}, see also \cite{pietraho:components}.
The work of  C.~Bonnaf{\'e} and L.~Iancu in the asymptotic parameter case relies on finding a generating set for cells defined in terms of standard bitableaux.
Finally, very recently  M.~Taskin has independently found another set of generators in the arbitrary rank case, see \cite{taskin:plactic}.

This paper is organized as follows.  In Section 2, we define the necessary objects and catalogue basic results.  Section 3 defines the family of operators $\Lambda^{r+1}$ and Section 4, in addition to showing that reducible combinatorial cells are stable under their action, describes the action on pairs of domino tableaux.   In Section 5, we verify the main result, and leave the proof of a few crucial lemmas to Sections 6 and 7.

\section{Definitions and Preliminaries}

\subsection{Robinson-Schensted Algorithms}

The hyperoctahedral group $W_n$ of rank $n$ is the group of permutations of
the set $\{\pm 1, \pm 2, \ldots , \pm n \}$ which commute with the
involution $i \mapsto -i$.  It is the Weyl group of type $B_n$.
We will write $w \in W_n$  in one-line notation as
$$w= (w(1)\; w(2) \ldots w(n)).$$

A Young diagram is a finite left-justified arrays of
squares arranged with non-increasing row lengths.  We will denote the square in row
$i$ and column $j$ of the diagram by  $S_{i,j}$ so that
$S_{1,1}$ is  the uppermost left square in the Young diagram below:
$$
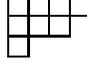
\begin{tiny}
\begin{tableau}
    :.{}.{}.{}.{}\\
    :.{}.{}.{}\\
    :.{}\\
\end{tableau}
\end{tiny}
$$
\begin{definition}
    Let $r \in \mathbb{N}$ and $\lambda$ be a partition of a positive integer $m$.
     A domino tableau of rank $r$ and shape $\lambda$ is a Young diagram
    of shape $\lambda$ whose squares are labeled by integers in such a way that $0$ labels $S_{ij}$ iff $i+j<r+2$, each element of some set $M$ labels exactly two adjacent squares, and all labels increase weakly along both rows and columns. A domino tableau is standard iff $M=\{1, \ldots, n\}$ for some $n$. We will write $SDT_r(n)$ for the set of standard domino tableaux with $n$ dominos.  The set of squares labeled by $0$ will be called the core of $T$.
\end{definition}

Following \cite{garfinkle1} and \cite{vanleeuwen:rank}, we
describe the Robinson-Schensted bijections
$$G_r : W_n \rightarrow SDT_r(n) \times SDT_r(n) $$ between elements of $W_n$ and same-shape pairs of rank $r$ standard domino tableaux. The algorithm
is based on an insertion map $\alpha$ which, given an entry $w(j)$ in the one-line notation for $w \in W_n$, inserts a domino with label $w(j)$ into a domino tableau.
This insertion map is
similar to the usual Robinson-Schensted insertion map and is
precisely defined in \cite{garfinkle1}(1.2.5). To construct the left
tableau, start with $T_1(0)$, the only tableau in $SDT_r(0)$. Define
$T_1(1) = \alpha (w(1),T_1(0))$ and continue
inductively by letting
$$T_1(k+1)=\alpha\big(w(k+1),T_1(k)\big).$$
The left domino tableau $T_1(n)$ will be standard and of rank $r$.
The right tableau keeps track of the sequence of shapes of the left tableaux; we define $T_2(n)$ to be the unique tableau so that $T_2(k) \in SDT_r(k)$ has the same shape as $T_1(k)$ for all $k \leq n$.  The Robinson-Schensted map is then defined by $G_r(w)=(T_1(n),T_2(n)).$ We will also often use the notation $(T_1(n-1),T_2(n-1)) = (T_1(n),T_2(n))'.$

\subsection{Cycles and Extended Cycles}  We briefly recall the definition of a cycle in a domino tableau as well as a number of related notions which we will use later.

For a standard domino tableau $T$ of arbitrary rank $r$, we will say the square $S_{ij}$ is fixed when $i+j$ has the opposite parity as $r$, otherwise,
we will call it variable.  If $S_{ij}$ is variable and $i$ is odd, we will say $S_{ij}$ is of type $X$; if $i$ is even, we will say $S_{ij}$ is of type $W$. We will write $D(k, T )$ for the domino labeled by the positive integer $k$ in $ T$ and $supp \,  D(k, T )$ will denote its underlying squares. Write $label \, S_{ij}$ for the label of the square $S_{ij}$ in $T$ . We extend
this notion slightly by letting $label \, S_{ij}$ = 0 if either $i$ or $j$ is less than or equal to
zero, and $label \, S_{ij} = \infty$ if $i$ and $j$ are positive but $S_{ij}$ is not a square in $T$ .

\begin{definition}
Suppose that  $supp \, D(k,T)= \{S_{ij},S_{i+1,j}\}$ or
$\{S_{i,j-1},S_{ij}\}$ and the square  $S_{ij}$ is fixed. Define
$D'(k)$ to be a domino labeled by the integer $k$ with $supp \,
D'(k,T)$ equal to
    \begin{enumerate}
        \item $\{S_{ij}, S_{i-1,j}\}$     if $k< label \, S_{i-1,j+1}$
        \item $\{S_{ij}, S_{i,j+1}\}$    if $k> label \, S_{i-1,j+1}$
    \end{enumerate}
Alternately, suppose that  $supp \, D(k,T)= \{S_{ij},S_{i-1,j}\}$
or $\{S_{i,j+1},S_{ij}\}$ and the square  $S_{ij}$ is fixed.
Define $supp \, D'(k,T)$ to be
    \begin{enumerate}
        \item $\{S_{ij},S_{i,j-1}\} $    if $k< label \, S_{i+1,j-1}$
        \item $\{S_{ij},S_{i+1,j}\}$         if $k> label \, S_{i+1,j-1}$
    \end{enumerate}
\end{definition}

\begin{definition}
The cycle $c=c(k,T)$ through $k$ in a standard domino tableau $T$ is
a union of labels of $T$  defined by the condition that $l \in c$ if
either
    \begin{enumerate}
        \item $l=k$,
        \item $supp \, D(l,T) \cap supp \, D'(m,T) \neq
        \emptyset$ for some $m \in c$, or
        \item $supp \, D'(l,T) \cap supp \, D(m,T) \neq
        \emptyset$ for some $m \in c$.
    \end{enumerate}
\end{definition}

Often, we will identify the labels contained in a cycle with their underlying dominos.
We define a domino tableau $MT(T,c)$ by replacing every domino $D(l,T) \in c$ by the corresponding domino $D'(l,T)$.  By \cite{garfinkle1}(1.5.27), the result is a standard tableau whose shape is either equal the
shape of $T$, or one square will be removed or added to the core,
and one will be added.
The cycle $c$ is called closed in the former
case and open in the latter.  If moving through $c$ adds a square to the core, we will call $c$ a core open cycle; the other open cycles will be called non-core.  For an open cycle $c$ of a tableau $T$, we will write $S_b(c)$ for the square that has been removed or
added to the core by moving through $c$.  Similarly, we will write
$S_f(c)$ for the square that is added to the shape of $T$. Note that
$S_b(c)$ and $S_f(c)$ are always variable squares.

\begin{definition}
A variable square $S_{ij}$ satisfying the conditions that
    \begin{enumerate}
        \item neither $S_{i,j+1}$ nor $S_{i+1,j}$ lie in $T$, and
        \item either
            \begin{enumerate}
                \item both $S_{i-1,j}$ and $S_{i,j-1}$ lie in $T$, or
                \item either $S_{i-1,j}$ lies in $T$ and $j=1$ or $S_{i,j-1}$ lies in $T$ and $i=1$,
            \end{enumerate}
    \end{enumerate}
will be called a hole if it is of type $W$ and a  corner if it is
of type $X$.  It will be called full if $S_{ij} \in T$ and empty otherwise.
\end{definition}

Let $U$ be a set of cycles in $T$. Because the order in which one moves through a set of cycles is immaterial by \cite{garfinkle1}(1.5.29), we can unambiguously write $MT(T,U)$ for the
tableau obtained by moving through all of the cycles in the set $U$.

Moving through a cycle in a pair of same-shape tableaux is somewhat problematic, as it may result in a pair of tableaux which in not same-shape.  We require the following
definition.

\begin{definition}
Consider $(T_1,T_2)$ a pair of same-shape domino tableaux, $k$ a label
in $T_1$, and $c$ the cycle in $T_1$ through $k$.  The extended cycle
$\tilde{c}$ of $k$ in $T_1$ relative to $T_2$ is a union of cycles in
$T_1$ which contains $c$.  Further, the union of two cycles $c_1 \cup
c_2$ lies in  $\tilde{c}$ if either one is contained in $\tilde{c}$ and,
for some cycle $d$ in $T_2$,  $S_b(d)$ coincides with a square of
$c_1$ and $S_f(d)$ coincides with a square of $MT(T_1,c_2)$.  The
symmetric notion of an extended cycle in $T_2$ relative to $T_1$ is
defined in the natural way.
\end{definition}

For an extended cycle $\tilde{c}$ in $T_2$ relative to $T_1$, write $\tilde{c} =c_1 \cup
\ldots \cup c_m$ and let  $d_1, \ldots, d_m$ be cycles in $T_1$ such that
$S_b(c_i)=S_b(d_i)$ for all $i$, $S_f(d_m)=S_f(c_1)$, and
$S_f(d_i)=S_f(c_{i+1})$ for $1\leq i < m$. The union $\tilde{d}=d_1
\cup \dots \cup d_m$ is an extended cycle in $T_1$ relative to $T_2$
called the extended cycle corresponding to $\tilde{c}$.
Symmetrically, $\tilde{c}$ is the extended cycle corresponding to
$\tilde{d}$.

We can now define a moving through operation for a pair of same-shape domino
tableaux.  Write $b$ for an ordered pair $(\tilde{c},\tilde{d})$
of extended cycles in $(T_1,T_2)$ that correspond to each other.  Define $MT((S,T), b)$ to equal  $(MT(S,\tilde{c}),MT(T,\tilde{d})).$  It is clear that
this operation produces another pair of same-shape domino
tableaux.  Speaking loosely, we will often refer to this operation as moving through either of the extended cycles $\tilde{c}$ or $\tilde{d}.$  Note that if the cycle $c$ is closed, then $\tilde{c}=c$ and moving through a pair of tableaux boils down to the operation $(MT(S,c),T).$

\begin{definition} We will say that a set of squares in a domino tableau is {\it boxed} iff it is entirely contained in a set of squares of the form $\{S_{ij},S_{i+1,j},S_{i,j+1},S_{i+1,j+1}\}$ where  $S_{ij}$ must be of type $X$.  A set will be called unboxed if it is not boxed.
\end{definition}

Boxing is well-behaved with respect to cycles.  If  $c$ is a cycle in $T$, then all of its underlying dominos are either boxed or unboxed.  Furthermore, moving through a domino changes its boxing, and consequently, the boxing of all the the dominos in the  cycle containing it.  The same holds for all of the dominos in an extended cycle of a same-shape tableau pair.

We will say that two sets of squares are adjacent in a tableau $T$ if there are two squares, one in each set, which share a common side.   One set of squares in $T$ will be said to be below another if all of its squares lie in rows strictly below the rows of the squares of the other set.  The notion of above is defined similarly.  We will often call the squares underlying a domino its position, and say that a position is extremal if its removal from a Young diagram results in another Young diagram corresponding to a partition of the same rank.

\subsection{Combinatorial Cells}

\begin{definition} Consider $w,v \in W_n$ of type $B_n$ and fix a non-negative integer
$r$ letting $G_r(w)=(T_1,T_2)$ and $G_r(v)=(\tilde{T}_1,\tilde{T}_2)$.  We will say that $w$ and $v$ are
    \begin{enumerate}
        \item in the same irreducible combinatorial left cell of rank $r$ if $T_2=\tilde{T}_2$, and
        \item in the same reducible combinatorial left cell of rank $r$ if
    there is a set $U$ of non-core open cycles in $T_2$ such that $\tilde{T}_2 = MT(T_2, U)$.
\end{enumerate}
 We will say that $w$ and $v$ are in the same irreducible and reducible combinatorial right cells iff their inverses lie in the same irreducible and reducible combinatorial left cells, respectively.
\end{definition}

When $r\geq n-1$, the situation is somewhat
simpler; there are no non-core open cycles implying that
combinatorial left cells are determined simply by right tableaux.  Furthermore,
by the main result of \cite{pietraho:rscore}, for these values of $r$ all
combinatorial cells are actually independent of $r$.

Reducible combinatorial left cells of rank $r$ can also be described in terms of right tableaux of ranks $r$ and $r+1$.  A similar characterization holds
for combinatorial right cells.
\begin{theorem}[\cite{pietraho:equivalence}]  Reducible combinatorial left
cells of rank $r$ in the Weyl group of type $B_n$  are generated by the
equivalence relations of having the same right tableau in either
rank $r$ or rank $r+1$.
\end{theorem}

In what follows we will focus on reducible combinatorial right cells. Since our generalizations of the Robinson-Schensted algorithm $G_r$ behave well with respect to inverses by merely changing the order of the two tableaux, all the statements can be easily modified to apply to reducible combinatorial left cells as well.

\section{Knuth Relations}
\label{section:knuth}
The purpose of this section is to define the operators on the hyperoctahedral group $W_n$ which generate the reducible combinatorial right cells of rank $r$.  We first recount the situation in type $A$ for the symmetric group $S_n$.

\subsection{Type A}

Writing the elements of $S_n$ in one-line notation,
the Knuth relations on $S_n$  are the transformations $$(w(1) \; w(2) \ldots w(j-2)\; w(j)\; w(j-1) \ldots w(n))$$ which transpose the $j$th and $(j+1)$st entries of $w \in S_n$ whenever
    \begin{enumerate}
        \item[(i)]
 $j \geq 2$ and $w(j-1)$ lies between $w(j)$ and $w(j+1)$, or
        \item[(ii)]  $j<n-1$ and $w(j+2)$ lies between $w(j)$ and $w(j+1)$.
    \end{enumerate}
Every Knuth relation preserves the Robinson-Schensted left tableau $T_1(w)$ of $w$.  However, even more can be said:

\begin{theorem}[\cite{knuth:art3}] Knuth relations generate the combinatorial left cells in $S_n$.  More precisely, given $w$ and $v$ with $T_1(w)=T_1(v)$, there is a sequence of Knuth relations sending $w$ to $v$.
\end{theorem}

Aiming to adapt this theorem to the hyperoctahedral groups, we begin by rephrasing the Knuth relations first in terms of the length function on $S_n$, and then again in terms of the $\tau$-invariant.

The group $S_n$ is a Weyl group of a complex semisimple Lie algebra $\mathfrak{g}$ of type $A_{n-1}$ with Cartan subalgebra $\mathfrak{h}$.
Let $\Pi_n =\{\alpha_1, \alpha_2, \ldots, \alpha_{n-1}\}$  be a set of simple roots for a choice of positive roots in  the root system $\Delta(\mathfrak{g},\mathfrak{h})$ and write $s_\alpha$ for the simple reflection corresponding to $\alpha \in \Pi_n$.
We view $S_n$ as the group generated by the above simple reflections and let $\ell: S_n \rightarrow \mathbb{Z}$ be the corresponding length function on $S_n$.  If we identify the reflection $s_{\alpha_j}$ with the transposition interchanging the $j$th and $(j+1)$st entries of the permutation corresponding to $w \in S_n$, then the Knuth relations on $S_n$  are the transformations taking $w$ to $ws_{\alpha_j}.$   Noting that $\ell(ws_{\alpha_j}) < \ell(w)$ exactly when $w(j) > w(j+1)$, the domain for the Knuth relations is the set of $w \in S_n$ satisfying
    \begin{enumerate}
        \item[(i)] $\ell(ws_{\alpha_i})<\ell(w)<\ell(ws_{\alpha_j})<\ell(ws_{\alpha_j}s_{\alpha_i})$, or
        \item[(ii)] $\ell(ws_{\alpha_i})>\ell(w)>\ell(ws_{\alpha_j})>\ell(ws_{\alpha_j}s_{\alpha_i})$
    \end{enumerate}
for some $\alpha_i \in \Pi_n$.   This condition is satisfied only when $\alpha_i$ and $\alpha_j$ are adjacent simple roots, that is, when $i=j \pm 1$.

Armed with this restatement of the Knuth conditions, following \cite{vogan:tau} we define the $\tau$-invariant for $S_n$ and a related family of operators.

\begin{definition} Write $W=S_n$.  For $w\in W$, let $\tau(w) = \{ \alpha \in \Pi_n \; | \; \ell(ws_\alpha) < \ell(w)\}$.  Given simple roots $\alpha$ and $\beta$ in $\Pi_n$, let
$$D_{\alpha \beta}(W) = \{w \in W \; | \; \alpha \notin \tau(w) \text{ and } \beta \in \tau(w)\}.$$ When $\alpha$ and $\beta$ are adjacent roots, define
$T_{\alpha \beta} : D_{\alpha \beta}(W) \rightarrow D_{\beta \alpha}(W)$ by
$$T_{\alpha \beta}(w) = \{ws_\alpha, ws_\beta\} \cap D_{\beta \alpha}(W).$$
\end{definition}
When defined, the operators $T_{\alpha \beta}$ are single-valued and preserve the Robinson-Schensted left tableau $T_1(w)$ of $w$.  The following is a direct consequence of the result on Knuth relations:
\begin{corollary}
The family of operators $T_{\alpha \beta} : D_{\alpha \beta}(W) \rightarrow D_{\beta \alpha}(W)$ generates the combinatorial left cells in $S_n$.  More precisely, given $w$ and $v$ with $T_1(w)=T_1(v)$, there is a sequence of $T_{\alpha \beta}$ operators sending $w$ to $v$.
\end{corollary}

\subsection{Type B}
\label{section:operators}
In order to define similar relations for the hyperoctahedral groups,  we mimic the final construction in type $A$.  The group $W_n$ is a Weyl group of a complex semisimple Lie algebra $\mathfrak{g}$ of type $B_n$ with Cartan subalgebra $\mathfrak{h}$.  Let $\{\epsilon_1, \ldots \epsilon_n\}$ be a basis for $\mathfrak{h}^*$ such that if we define $\alpha'_1=\epsilon_1 $ and $\alpha_i = \epsilon_i-\epsilon_{i-1},$  then $$\Pi_n = \{\alpha'_1, \alpha_2, \ldots, \alpha_n\}$$ is the set of simple roots for some choice of positive roots $\Delta^+(\mathfrak{g},\mathfrak{h})$.  While this choice of simple roots is not standard, we adopt it to obtain somewhat cleaner statements and reconcile our work with \cite{garfinkle3}.

 We modify $\Pi_n$ slightly to include certain non-simple roots.  For $i \leq n$, let $\alpha'_i=\alpha'_1+\alpha_2 + \ldots \alpha_i$ and when $k$ is a non-negative integer write $$\Pi_n^k = \{\alpha'_1, \alpha'_2, \ldots \alpha'_{\min(n,k)}, \alpha_2, \ldots \alpha_n\}.$$
Further, write $s_i$ for the simple reflection $s_{\alpha_{i+1}}$ and $t_i$ for the reflection $s_{\alpha'_{i}}$.  We realize $W_n$ as a set of signed permutations on $n$ letters by identifying $s_i$ with the transposition $(1  \; 2 \ldots i+1 \; i \ldots n)$ and $t_i$ with the element $(1  \; 2 \ldots -i \ldots n).$

The generating set for reducible combinatorial right cells will be drawn out the following three types of operators:

\begin{definition}
For $w\in W_n$, and a non-negative integer $k$, let $$\tau^{k}(w) = \{ \alpha \in \Pi^{k}_n \; | \; \ell(ws_\alpha) < \ell(w)\}.$$  Given roots $\alpha$ and $\beta$ in $\Pi^{k}_n$ define
$D^{k}_{\alpha \beta}(W_n) = \{w \; | \; \alpha \notin \tau^{k}(w) \text{ and } \beta \in \tau^{k}(w)\}.$ The operator
$T^k_{\alpha \beta} : D^{k}_{\alpha \beta}(W_n) \rightarrow D^{k}_{\beta \alpha}(W_n)$ is defined by
$$T^{k}_{\alpha \beta}(w) = \{ws_\alpha, ws_\beta\} \cap D^{k}_{\beta \alpha}(W_n).$$
When $k=0$, we will eliminate the superscript and write $T_{\alpha \beta}$ for $T^0_{\alpha \beta}$ and $D_{\alpha \beta}(W_n)$ for $D^0_{\alpha \beta}(W_n).$
\end{definition}

\begin{definition}
We will say $w \in D^k_{I\!N}(W_n)$ if for some  $\{\delta,\beta\} = \{\alpha_{k+1},\alpha'_k\}$, we have
$$w \in  D^k_{\delta \beta}(W_n) \text{ and } ws_k \in D^k_{\beta \delta}(W_n).$$
Define $T^{k}_{I\!N}$ for  $w \in D^k_{I\!N}(W_n)$ by $T^{k}_{I\!N}(w) = w s_k.$
\end{definition}

\begin{definition}
We will say $w \in D^k_{SC}(W_n)$ if for some choice of
$\{\delta_i,\beta_i\} = \{\alpha_{i+1},\alpha'_i\}$ for every $i \leq k$, we have
$$w \in D^k_{\delta_i \beta_i}(W_n) \text{ and } wt_i \in D^k_{\beta_i \delta_i}(W_n) \hspace{.07in} \forall i \leq k.$$
Define $T^k_{SC}$ for  $w \in D^k_{SC}(W_n)$ by $T^{k}_{SC}(w) = w t_1.$
\end{definition}

To phrase the above in more digestible terms, we note that in $W_n$ the length function satisfies
 $\ell(ws_j)< \ell(w)$ iff $w(j)>w(j+1)$, and
 $\ell(wt_j)< \ell(w)$ iff $w(j)<0$.
Armed with this characterization, the domains and actions for the three types of operators which we will be interested in can be described more succinctly.
\begin{itemize}
    \item[$T_{\alpha \beta}:$] $w$ lies in $D_{\alpha \beta}(W)$ for
    $\{\alpha, \beta\}= \{\alpha_{j+1},\alpha_{j+2}\}$ whenever $w(j+1)$ is either greater than or smaller than both $w(j)$ and $w(j+2)$; $T_{\alpha \beta}$ then interchanges the smallest and largest entries among $w(j), w(j+1),$ and $w(j+2).$
        \item[$T^{k}_{I\!N}:$] $w \in D^{k}_{I\!N}(W)$ iff $w(k)$ and $w(k+1)$ are of opposite sign.  The operator $T^{k}_{I\!N}$ interchanges the $k$th and $k+1$st entries of $w$.
            \item[$T^{k}_{SC}:$] $w \in D^{k}_{SC}(W)$ iff $|w(1)| > |w(2)| > \ldots > |w(k+1)|$.  The operator $T^{k}_{SC}$ changes the sign of $w(1)$.
\end{itemize}

\begin{definition}  For an integer $k$, we define a set of operators
$$\Lambda^{k} = \{T_{\alpha \beta} \; | \; {\alpha, \beta} \in \Pi^{0}_n \} \cup \{T^{i}_{I\!N}  \; | \; i \leq k \} \cup \{T^{k}_{SC}\}.$$
\end{definition}

This is the sought-after set of generators for reducible combinatorial right cells of rank $r$.  The following are our two main results:

\begin{theorem} The operators in $\Lambda^{r+1}$ preserve reducible combinatorial left cells of rank $r$.  If $S \in \Lambda^{r+1}$, and $w \in W_n$ is in the domain of $S$,  then the Robinson-Schensted left domino tableaux  $T_1(w)$ and $T_1(Sw)$ of rank $r$  differ only by moving through a set of non-core open cycles in $T_1(w)$.
\label{theorem:stable}
\end{theorem}

\begin{theorem} The family of operators $\Lambda^{r+1}$ generates the reducible combinatorial left cells of rank $r$.  More precisely, given $w$ and $v \in W_n$ whose Robinson-Schensted left domino tableaux  $T_1(w)$ and $T_1(v)$ of rank $r$ differ only by moving through a set of non-core open cycles, there is a sequence of operators in $\Lambda^{r+1}$ sending $w$ to $v$.
\label{theorem:generators}
\end{theorem}

These results have previously been obtained in two special cases.  In her work on the primitive spectrum of the universal enveloping algebras in types $B$ and $C$, Devra Garfinkle constructed a generating set for the reducible combinatorial cells of rank zero \cite{garfinkle3}.  Additionally, C.~Bonnaf{\'e} and L.~Iancu obtained a generating set for combinatorial left cells in the so-called asymptotic case when $r \geq n-1$.  The generating set $\Lambda^{r+1}$ proposed above for domino tableaux of arbitrary rank generalizes both of these results.

\begin{theorem}[D.~Garfinkle]  The reducible combinatorial left cells of rank zero are generated by the family of operators $T^1_{\alpha \beta}$ where $\alpha$ and  $\beta$ are adjacent simple roots in  $ \Pi_n.$
\end{theorem}

At first glance, the set $\Lambda^1$ is not exactly  Garfinkle's generating set.   The latter contains multi-valued operators $T^1_{\alpha'_1 \alpha_2}$ and $T^1_{\alpha_2 \alpha'_1}$.  However, we note that the  combined action and domains of $T^1_{\alpha'_1 \alpha_2}$ and $T^1_{\alpha_2 \alpha'_1}$ agree precisely with those of $T^1_{SC}$ and $T^1_{I\!N} \in \Lambda^1.$ The operators $T^1_{SC}$ and $T^1_{I\!N}$ merely split up Garfinkle's original multi-valued operators.

\begin{theorem}[C.~Bonnaf{\'e} and L.~Iancu.]  When the rank $r \geq n-1$, the reducible and irreducible combinatorial left cells coincide.  They are generated by the family of operators $\Lambda^{r+1} = \Lambda^{n}.$
\label{theorem:bi}
\end{theorem}

We note that in the asymptotic case, the operator $T^{r+1}_{SC}$ is not defined, so that $\Lambda^{r+1}$ consists entirely of the $T_{\alpha \beta}$ and $T^k_{I\!N}$ operators, as exhibited in \cite{bonnafe:iancu}.

It was observed in \cite{pietraho:constructible} that combinatorial right cells are in general not well-behaved with increasing rank.  This phenomenon is readily explained by examining the composition of $\Lambda^{r+1}$ more closely.  While increasing $r$ expands the family of  $T_{I\!N}$ operators contained in $\Lambda^{r+1}$, it also diminishes the domain of $T_{SC}$.   The complicated behavior of right cells as $r$ increases is
 a manifestation of this interplay.

\section{Stability under $\Lambda^{r+1}$}

The goal of this section is to prove Theorem \ref{theorem:stable}, showing that reducible right cells in rank $r$ are invariant under the operators in $\Lambda^{r+1}$. We will rely on a decomposition of $W_n$, which we describe presently.
 Let $W_m$ be the subgroup of $W_n$ generated by the  reflections $\{s_1, \ldots, s_{m-1}, t_1\}$
and define $X_m^n$ to be the set of $x \in W_n$ which satisfy
$$0<x(1)<x(2)< \ldots < x(m).$$
Then $X_m^n$ is a cross-section of $W_n/W_m$ and we can write every $w \in W_n$ as a product $w=xw'$ with $x \in X_m^n$ and $w'\in W_m$.  The following appears in \cite{bgil}.

\begin{proposition} If $x \in X_m^n$ and $w', v' \in W_m$ are in the same irreducible combinatorial left cell of rank $r$, then  $xw'$ and $xv' \in W_n$ are in the same irreducible combinatorial left cell of rank $r$.
\label{proposition:cross}
\end{proposition}

\begin{corollary}
If $x \in X_m^n$ and $w', v' \in W_m$ are in the same reducible combinatorial left cell of rank $r$, then  $xw'$ and $xv' \in W_n$ are in the same reducible combinatorial left cell of rank $r$.
\label{corollary:cross}
\end{corollary}

\begin{proof}  The main result of \cite{pietraho:equivalence} shows that reducible combinatorial left cells of rank $r$ are the least common refinements of irreducible combinatorial left cells of ranks $r$ and $r+1$.  Hence the corollary follows directly from the above proposition.
\end{proof}

\subsection{The Operator $T_{SC}$}
We first examine the family of operators $T^k_{SC}$ which, under appropriate circumstances, change the sign of the first entry of $w$.

\begin{proposition} The operator $T^{r+1}_{SC}$ preserves  the irreducible and, consequently, the reducible combinatorial right cells of rank $r$ in $W_{r+2}$.
\label{proposition:SC}
\end{proposition}

\begin{proof}  Consider $w \in D^{r+1}_{SC}(W_{r+2})$ so that
$$|w(1)| > |w(2)| > \ldots |w(r+2)|.$$  Write $v=w^{-1}$ and $u=(T^{r+1}_{SC}(w))^{-1}.$
We will show that their right tableaux $T_2(v)$ and  $T_2(u)$ agree. Note that the condition on $w$ implies $|v(1)| > |v(2)| > \ldots |v(r+2)|$ and the fact that $T^{r+1}_{SC}$ only changes the sign of $w(1)$ forces
    \begin{enumerate}
        \item $v(i)=u(i)$ for $i\leq r+1$, and
        \item $u(r+2)=-v(r+2)\in \{\pm 1\}$.
    \end{enumerate}
If we write the positive entries in $v$ as $a_1 > \ldots > a_p$ and the negative ones as $-b_1< \ldots<  -b_q$, then the left tableau of $v$ {\it before} the insertion of $v(r+2)$ must have the form:

$$
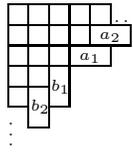
\begin{tiny}
\begin{tableau}
:.{}.{}.{}.{}.{} \cdots\\
:.{}.{}.{}.{}>{a_2}\\
:.{}.{}.{}>{a_1}\\
:.{}.{}^{\,b_1}\\
:.{}^{\,b_2}\\
: \vdots \\
\end{tableau}
\end{tiny}
$$
This is also the left tableau of $u$ before the insertion of $u(r+2)$ and we point out that the corresponding right tableaux for both $u$ and $v$ are the same.  If $a_1>b_1$, then the insertion of $v(r+2)$ and $u(r+2)$ into the above tableau yields either
$$
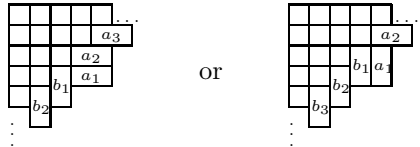
\begin{tiny}
\begin{tableau}
:.{}.{}.{}.{}.{} \cdots\\
:.{}.{}.{}.{}>{a_3}\\
:.{}.{}.{}>{a_2}\\
:.{}.{}^{\,b_1}>{a_1}\\
:.{}^{\,b_2}\\
: \vdots \\
\end{tableau}
\end{tiny}
\hspace{.2in}
\raisebox{.35in}{
\text{ or }}
\hspace{.3in}
\begin{tiny}
\begin{tableau}
:.{}.{}.{}.{}.{} \cdots\\
:.{}.{}.{}.{}>{a_2}\\
:.{}.{}.{}^{\,b_1}^{\,\,a_1}\\
:.{}.{}^{\,b_2}\\
:.{}^{\,b_3}\\
: \vdots \\
\end{tableau}
\end{tiny}
$$
In either case, the same domino is added to their shared right tableau and $T_2(v) = T_2(u)$.   The proof is similar when $a_1<b_1$.
\end{proof}

\begin{corollary}When $n \geq r+2$, the operator $T^{r+1}_{SC}$ preserves  the irreducible and, consequently, the reducible, combinatorial left cells of rank $k$ in $W_{n}$.
\end{corollary}
\begin{proof}
For $w \in D^{r+1}_{SC}(W_n)$, write $w=xw' \in X_n^{r+2} W_{r+2}$.  Then $w' \in D^{r+1}_{SC}(W_{r+2})$ and by Proposition \ref{proposition:SC}, $T^{r+1}_{SC}(w')$ and $w'$ share the same irreducible combinatorial left cell of rank $k$.  Since $T^{r+1}_{SC}(w)= x T^{r+1}_{SC}(w') \in X_n^{r+2} W_{r+2}$, Proposition \ref{proposition:cross} and Corollary \ref{corollary:cross} imply the same is true of $T^{r+1}_{SC}(w)$ and $w$.
\end{proof}

Our final aim is to describe the action of $T^{r+1}_{SC}$ on the right tableau of $w$.  This is more or less captured in the proof of the above proposition, but we will attempt to be more precise.  Consider $w \in D^{r+1}_{SC}(W_{n})$, write $T_2$ for the right tableau $T_2(w)$ of rank $r$, and define $T_2(r+2)$ to be the subtableau of $T_2$ which contains only dominos whose labels are less than or equal to $r+2$.  Because  $w \in D^{r+1}_{SC}(W_{n})$, $T_2(r+2)$ contains
exactly  the dominos adjacent to the core of $T_2$ as well as the unique domino of $T_2$ which shares its long edge with the long edge of one of the dominos adjacent to the core.  Furthermore, $T_2(r+2)$ must have the form
$$
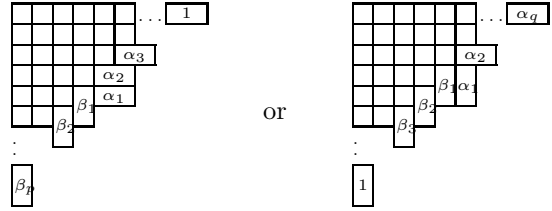
\begin{tiny}
\begin{tableau}
:.{}.{}.{}.{}.{}.{} \cdots >{1}\\
:.{}.{}.{}.{}.{}.{}\\
:.{}.{}.{}.{}.{}>{\alpha_3}\\
:.{}.{}.{}.{}>{\alpha_2}\\
:.{}.{}.{}^{\,\,\beta_1}>{\alpha_1}\\
:.{}.{}^{\,\,\beta_2}\\
:\vdots\\
:^{\,\,\beta_p}\\
\end{tableau}
\end{tiny}
\hspace{.2in}
\raisebox{.35in}{
\text{ or }}
\hspace{.3in}
\begin{tiny}
\begin{tableau}
:.{}.{}.{}.{}.{}.{} \cdots >{\alpha_q}\\
:.{}.{}.{}.{}.{}.{} \\
:.{}.{}.{}.{}.{}>{\alpha_2}\\
:.{}.{}.{}.{}^{\,\,\beta_1}^{\,\,\alpha_1}\\
:.{}.{}.{}^{\,\,\beta_2}\\
:.{}.{}^{\,\,\beta_3}\\
:\vdots\\
:^{1}\\
\end{tableau}
\end{tiny}
$$

\vspace{.2in}
\noindent
where $\alpha_1=r+2$, noting that $p$ may equal zero.  By the proof of the preceding proposition, $T^{r+1}_{SC}$ acts on $T_2(r+2)$ by interchanging one of its possible two forms of the same shape with the other.  It is easy to see that the remaining dominos of the tableaux $T_2$ and the right tableau $T_2(T^{r+1}_{SC}(w))$ are the same, since both keep track of the subsequent insertions into the left tableau $T_1(w)(r+2) = T_1(T^{k+1}_{SC}(w))(r+2).$

\begin{example}
Let $w=(4 \; -3 \; -2 \; 1) \in D^3_{SC}(W_4).$ Then $T^3_{SC}(w)=(-4 \; -3 \; -2 \; 1)$ and
$$
\raisebox{.1in}{$T_2(w)=$}
\hspace{.2in}
\begin{tiny}
\begin{tableau}
:.{}.{}>{1}\\
:.{}^{3}>4\\
:^2\\
\end{tableau}
\end{tiny}
\hspace{.2in}
\raisebox{.1in}{
\text{ and  \hspace{.1in} }  $T_2(T^3_{SC}(w))=$}
\hspace{.2in}
\begin{tiny}
\begin{tableau}
:.{}.{}^{3}^{4}\\
:.{}^{2}\\
:^{1}\\
\end{tableau}
\end{tiny}
\hspace{.1in}
\raisebox{.1in}{.}
$$

\end{example}

On a final note, we observe that this description of the action of $T^{r+1}_{SC}$ reproduces the action of a portion of the multi-valued $T_{\alpha'_1 \alpha_2}$ operator detailed in \cite{garfinkle2}(2.3.4) for the rank zero case.

\subsection{The Operators $T_{\alpha \beta}$}
\label{section:alphabeta}
Given adjacent simple roots $\{\alpha, \beta\}=\{\alpha_{j+1},\alpha_{j+2}\}$, an operator $T_{\alpha \beta}$  interchanges the smallest and largest entries among $w(j)$, $ w(j+1),$ and $w(j+2).$  The following is a consequence of Proposition 3.10 in \cite{bgil}.

\begin{proposition} The operators $T_{\alpha \beta}$ preserve  the irreducible and, consequently, the reducible combinatorial right cells of rank $r$ in $W_{n}$.
\label{proposition:alphabeta}
\end{proposition}

What is still required is a description of their action on right tableaux.  We begin with an explicit description of their domains.  Recall that $\alpha_{j+1} \in \tau(w)$ iff $\ell(ws_{\alpha_{j+1}}) < \ell(w)$ which occurs iff $w(j)>w(j+1)$.  This condition is easily read off from the right tableau $T_2(w)$.  We will say that $k$ lies below $l$ in a domino tableau iff every row containing a square of the domino labeled $k$ lies below  every row of the domino labeled $l$.  Then $\alpha_{j+1} \in \tau(w)$ iff $j+1$ lies below $j$ in $T_2(w)$.  Unraveling the definition of $D_{\alpha \beta}(W)$, we find that when $\{\alpha, \beta\}=\{\alpha_{j+1},\alpha_{j+2}\}$, $w$ lies in the domain of $T_{\alpha \beta}$ iff in the tableau $T_R(w)$ either:
    \begin{enumerate}
        \item $j+1$ lies below $j$ and $j+2$ does not lie below $j+1$, or
        \item $j+1$ does not lie below $j$ and $j+2$ lies below $j+1$.
    \end{enumerate}
Two cases are necessary to describe the action of $T_{\alpha \beta}$.

\vspace{.1in}
{\it Case 1.}\hspace{.05in}
In the following, let $k=j+1$ and $l=j+2$.  If either one of the following four configurations of dominos appears in $T_2(w)$:
$$
\raisebox{.1in}{$F_1(j)=$}
\hspace{.2in}
\begin{tiny}
\begin{tableau}
:^{j}>{k}\\
:;>{l}\\
\end{tableau}
\end{tiny}
\hspace{.2in}
\raisebox{.1in}{
\text{ or  \hspace{.1in} }  $\widetilde{F}_1(j)=$}
\hspace{.2in}
\begin{tiny}
\begin{tableau}
:>{j}^{l}\\
:>{k}\\
\end{tableau}
\end{tiny}
\hspace{.1in}
\raisebox{.1in}{,}
$$

$$
\raisebox{.1in}{$F_2(j)=$}
\hspace{.2in}
\begin{tiny}
\begin{tableau}
:>{j}\\
:^k^{l}\\
:;\\
\end{tableau}
\end{tiny}
\hspace{.2in}
\raisebox{.1in}{
\text{ or  \hspace{.1in} }  $\widetilde{F}_2(j)=$}
\hspace{.2in}
\begin{tiny}
\begin{tableau}
:^{j}^{k}\\
:;\\
:>{l}\\
\end{tableau}
\end{tiny}
\hspace{.1in}
\raisebox{.1in}{,}
$$
then the action $T_{\alpha \beta}$ on the right tableau of $w$ swaps $F_i(j)$ and $\widetilde{F}_i(j)$ within $T_R(w)$.

  \vspace{.1in}
{\it Case 2.}\hspace{.05in}
If none of the above four configurations appear in $T_2(w)$ and $w$ lies in the domain $D_{\alpha \beta}(W_n)$, then $T_{\alpha \beta}$ acts by swapping either the dominos labeled $j$ and $j+1$ or $j+1$ and $j+2$.

The proof that this description of $T_{\alpha \beta}$ on domino tableau accurately depicts the action of $T_{\alpha \beta}$ defined on $W_n$ is not difficult.  It appears as \cite{garfinkle2}(2.1.19) in rank zero, and follows almost identically for higher rank tableaux.

\subsection{The Operators $T^k_{I\!N}$}

The operator $T^k_{I\!N}$  interchanges the entries $w(k)$ and $w(k+1)$ in $w$ whenever they are of opposite signs.  The following is again a consequence of Proposition 3.10 in \cite{bgil}.

\begin{proposition} When $k \leq r$, the operator $T^k_{I\!N}$ preserves  the irreducible and, consequently, the reducible combinatorial left cells of rank $r$ in $W_{n}$.
\label{proposition:in}
\end{proposition}

Slightly more is true:

\begin{proposition} The operator $T^{r+1}_{I\!N}$ preserves the reducible combinatorial left cells of rank $r$ in $W_{n}$.
\end{proposition}
\begin{proof}By Proposition \ref{proposition:in}, the operator $T^{r+1}_{I\!N}$ preserves the irreducible combinatorial left cells of rank $r+1$ in $W_{n}$.  Since reducible combinatorial left cells of rank $r$ are the least common refinements of irreducible cells in ranks $r$ and $r+1$, the result follows.
\end{proof}

Next, we describe the action of $T^k_{I\!N}$ on the left and right tableaux  $T_1(w)$ and $T_2(w)$.
Recall that $\alpha'_k \in \tau(w)$ iff $\ell(wt_k) < \ell(w)$ which occurs iff $w(k)<0$.  As long as $k \leq r+1$, this occurs iff the domino with label $k$ in $T_2(w)$ is vertical.  Recalling that $\alpha_{k+1} \in \tau(w)$ iff $k+1$ lies below $k$ in $T_2(w)$, we can describe the domain of $T^k_{I\!N}$ as follows.

\begin{definition} We will call a tableau $T$ of rank $r$ {\it sparse} if there is a square of the form $S_{m,r+3-m}$ which is empty in $T$.
\end{definition}

\vspace{.1in}
{\it Case 1.}\hspace{.05in}
If the tableau $T_2(w)(r+2)$ is sparse, then  $w \in S^k_{I\!N}$ for all $k\leq r+1$ iff  the dominos with labels $k$ and $k+1$ in $T_2(w)$ are of opposite orientations.  In this setting,  since $T^k_{I\!N}$ swaps $w(k)$ and $w(k+1)$ in $w$ whenever they are of opposite sign, its action merely reverses the order in which the $k$th and $k+1$st dominos are inserted into the left tableau.  Since $T_2(w)(r+2)$ is sparse, these insertions do not interact with each other. Thus on the shape-tracking tableau $T_2(w)$, $T^k_{I\!N}$ merely swaps the dominos with labels $k$ and $k+1$, while acting trivially on $T_1(w)$.

\vspace{.1in}
{\it Case 2.}\hspace{.05in}
If the tableau $T_2(w)(r+2)$ is not sparse, then $k=r+1$
and the insertions of $w(k)$ and $w(k+1)$ into the left tableau interact with each other.

It is easy to see that there are four possible configurations of dominos with labels $k=r+1$ and $l=r+2$ within $T_2(w)$ when $w \in S^{r+1}_{I\!N}(W)$:

$$
\raisebox{.1in}{$E_0(k)=$}
\hspace{.2in}
\begin{tiny}
\begin{tableau}
:^k^l\\
:;\\
\end{tableau}
\end{tiny}
\hspace{.2in}
\raisebox{.1in}{
\text{ or  \hspace{.1in} }  $\widetilde{E}_0(k)=$}
\hspace{.2in}
\begin{tiny}
\begin{tableau}
:^{k}>l\\
:;\\
\end{tableau}
\end{tiny}
\hspace{.1in}
\raisebox{.1in}{,}
$$

$$
\raisebox{.1in}{$E_1(k)=$}
\hspace{.2in}
\begin{tiny}
\begin{tableau}
:>k\\
:>l\\
\end{tableau}
\end{tiny}
\hspace{.2in}
\raisebox{.1in}{
\text{ or  \hspace{.1in} }  $\widetilde{E}_1(k)=$}
\hspace{.2in}
\begin{tiny}
\begin{tableau}
:>k\\
:^l\\
:;\\
\end{tableau}
\end{tiny}
\hspace{.1in}
\raisebox{.1in}{.}
$$

\begin{proposition}
Consider $w \in D^{r+1}_{I\!N}(W)$ and let $v= T^{r+1}_{I\!N}(w)$.  If
$G_r(w)=(T_1, T_2)$, then the tableau pair  $G_r(v)=(\tilde{T}_1, \tilde{T}_2)$ admits the following description:

\begin{enumerate}
    \item If $T_2(r+2)$ is sparse, then $T_1=\tilde{T}_1$ and $\tilde{T}_2$ is obtained by swapping the dominos  with labels $r+1$ and $r+2$ in the tableau $T_2$.

    \item If $T_2(r+2)$ is not sparse, then the description of the action depends on the exact configuration of the dominos with labels $r+1$ and $r+2$ in $T_2$:
    \begin{enumerate}
        \item If the configuration $E_0(r+1)$ or $E_1(r+1)$ appears in $T_2$, write it as $E_i(r+1)$ and  let $\bar{T}_2 = (T_2 \setminus E_i(r+1))\cup E_{1-i}(r+1)$.  Let $c$ be the extended cycle through $r+2$ in $\bar{T}_2$ relative to $T_1$.  Then
            $$(\tilde{T}_1, \tilde{T}_2)  = MT((T_1,\bar{T}_2), c).$$
        \item If $\widetilde{E}_0(r+1)$ or $\widetilde{E}_1(r+1)$ appears in $T_2$, let $c$ be the extended cycle through $r+2$ in $T_2$ relative to $T_1$ and define
                $(\bar{T}_1, \bar{T}_2) = MT((T_1,T_2), c).$
            Note that $\bar{T}_2$ must contain one of the configurations $E_0(k)$ or $E_1(k)$, which we label $E_i(r+1)$.  Then
            $$ (\tilde{T}_1, \tilde{T}_2) = \big(\bar{T}_1, (\bar{T}_2 \setminus E_i(r+1))\cup E_{1-i}(r+1)\big).$$
    \end{enumerate}
\end{enumerate}
\label{proposition:IN}
\end{proposition}

\begin{proof}
The first case has already been considered above.
In the special situation when $n=r+2$, the second case follows by inspection.  We mimic the proof of \cite{garfinkle2}(2.3.8) to verify the second case in general.
So consider $w \in D^{r+1}_{I\!N}(W)$
assuming that $T_2(r+2)$ is not sparse.
By symmetry and the fact that $T^{r+1}_{I\!N}$ is an involution, it is sufficient to consider only $w$ for which also $|w(r+1)| < |w(r+2)|$, $w(r+1)>0$ and $w(r+2)<0$. For such a $w$, let $\bar{w}=vt_{r+1}$.  Write $G_r(\bar{w})=(\bar{T}_1,\bar{T}_2).$
Then
    \begin{enumerate}
        \item $\widetilde{E}_1({r+1}) \subset T_2$, $E_0({r+1}) \subset \tilde{T}_2$, and $E_1(r+1) \subset \bar{T}_2,$
        \item $\bar{T}_2 = (\tilde{T}_2\setminus E_0({r+1}))\cup E_1({r+1})$ and $\bar{T}_1 = \tilde{T}_1,$
        \item $(\bar{T}_1,\bar{T}_2) = MT((T_1,T_2), c)$ where $c$ is the extended cycle through $r+2$ in $T_2$ relative to $T_1$.
    \end{enumerate}
Once these are verified, the proposition follows.  The last two parts imply that  $\tilde{T}_1=\bar{T}_1$, and the latter equals $MT(T_1,d)$ where $d$ is the extended cycle in $T_1$ corresponding to $c$, as desired.  Meanwhile, the second part then implies that $\tilde{T}_2$ will be as specified by the proposition.

Statements (1) and (2) follow easily from the definition of domino insertion, while the proof of (3) is identical to the rank zero case.  That the extended cycle through $r+2$ in $T_2$  consists only of non-core cycles follows from the remark at the end of this section, and then, for non-core cycles, the description detailing the relationship between moving through and domino insertion of \cite{garfinkle2}(2.3.2) carries without modification to the arbitrary rank case.

\end{proof}

\begin{example}
Let $w=(4 , -3 , -2 , 1) \in D^3_{I\!N}(W_4).$ Then $T^3_{I\!N}(w)=(4, -3, 1, -2)$ and the corresponding tableau pairs are:
$$
\raisebox{.1in}{$G_2(w)=$}
\hspace{.1in}
\raisebox{.05in}{\Bigg(}
\hspace{.05in}
\begin{tiny}
\begin{tableau}
:.{}.{}>{1}\\
:.{}^{3}>4\\
:^2\\
\end{tableau}
\end{tiny}
,
\begin{tiny}
\begin{tableau}
:.{}.{}>{1}\\
:.{}^{3}>4\\
:^2\\
\end{tableau}
\end{tiny}
\hspace{.05in}
\raisebox{.05in}{\Bigg)}
\hspace{.1in}
\raisebox{.1in}{
\text{ and  }  $G_2(T^3_{I\!N}(w))=$}
\hspace{.1in}
\raisebox{.05in}{\Bigg(}
\hspace{.05in}
\begin{tiny}
\begin{tableau}
:.{}.{}>{1}\\
:.{}^{3}^4\\
:^2\\
\end{tableau}
\end{tiny}
,
\begin{tiny}
\begin{tableau}
:.{}.{}>{1}\\
:.{}>{3}\\
:^2>4\\
\end{tableau}
\end{tiny}
\hspace{.05in}
\raisebox{.05in}{\Bigg)}
\raisebox{.1in}{.}
$$
The extended cycle through $r+2=4$ in $T_2(w)$ consists of the open cycle $\{4\}$ in $T_2(w)$ and the corresponding open cycle $\{4\}$ in $T_1(w).$
\end{example}

\begin{remark} In the case when $T_2(r+2)$, is not sparse it is not immediately clear that the operation on tableaux described in Proposition \ref{proposition:IN} actually produces a domino tableau of rank $r$. That it does follows from the fact that the extended cycle through which the tableaux are being moved through does not contain any core open cycles.  We verify this presently.  Suppose that $w \in T^{r+1}_{I\!N},$ so that $w(r+1)$ and $w(r+2)$ are of opposite sign.

{\it Case 1.}\hspace{.05in} Suppose either $\widetilde{E}_0(r+1)$ or     $\widetilde{E}_1(r+1)$ appears in $T_2$, and without loss of generality, assume that it is $\widetilde{E}_0(r+1)$.  The domino with label $r+2$ is adjacent to $r+1$ as well as to another domino adjacent to the core of $T_2$, which we label $l$.  For $r+2$ to be in an extended cycle containing a core open cycle, the extended cycle must also contain either $r+1$ or $l$.  We show that this is impossible.
First, $r+1$ and $r+2$ cannot be in the same extended cycle;  one of them is boxed but not the other. Since $w(r+1)$ and $w(r+2)$ are of opposite sign, it is easy to check that {\it within} $T_2(r+2)$, $r+2$ and $l$ are in different extended cycles.  But Lemmas 3.6 and 3.7 of \cite{pietraho:rscore} imply that this remains the case  within the full tableau $T_2$, see also \cite{garfinkle2}(2.3.3).

{\it Case 2.}\hspace{.05in}Suppose either ${E}_0(r+1)$ or     ${E}_1(r+1)$ appears in $T_2$, and without loss of generality, assume that it is ${E}_0(r+1)$.  The first operation prescribed by Proposition \ref{proposition:IN} is to swap ${E}_0(r+1)$ with
  ${E}_1(r+1)$ in $T_2$.  If we now define $l$ as before, then $r+2$ and $l$ must be in different extended cycles as one of them is boxed but not the other.  Since $w(r+1)$ and $w(r+2)$ are of opposite sign, it is again easy to check that {\it within} $(T_2(r+2)\setminus{E}_0(r+1)) \cup {E}_1(r+1) $, $r+2$ and $r+1$ are in different extended cycles, and the rest of the proof follows as above.
\label{remark:ec}
\end{remark}

\section{Generators of Combinatorial Cells}

The goal of this section is to verify that the set of operators in $\Lambda^{r+1}$ generates the reducible combinatorial right cells of rank $r$ in the hyperoctahedral group $W_n.$    We state Theorem \ref{theorem:generators} more precisely as:

\begin{theorem}  Suppose that $w,v \in W_n$ with $G_r(w)=(T_1, T_2)$ and $G_r(v) = (\tilde{T}_1, \tilde{T}_2)$.  If $\tilde{T}_1=MT(T_1,U)$ for some set of non-core open cycles $U$ in $T_1$, then there is a sequence of operators $\Sigma$ in the set $\Lambda^{r+1}$ such that $\Sigma(w)=v.$
\end{theorem}

Its proof follows directly from two auxiliary facts.  However, we need a definition first.

\begin{definition} We will say that a domino tableau is {\it somewhat special} if all of its non-core open cycles are boxed.
\end{definition}

A tableau is somewhat special iff all of its corners are empty, thus this is really a property of the underlying partition. Since moving through a non-core open cycle changes it from boxed to unboxed and vice-versa, there is a unique somewhat special tableau in the family of tableaux obtained from each other by moving through non-core open cycles.  We will write $S(T)$ for the somewhat special tableau corresponding to $T$.

\begin{proposition} Given $w\in W_n$, there is a sequence of operators $\Sigma$ in  $\Lambda^{r+1}$ such that the left tableau of $\Sigma(w)$ is somewhat special.
\label{prop:makespecial}
\end{proposition}

\begin{proposition} Suppose that $w,v \in W_n$ with $G_r(w)=(T_1, T_2)$ and $G_r(v) = (\tilde{T}_1, \tilde{T}_2)$.  If $T_1 = \tilde{T}_1$ and both are somewhat special, then there is a sequence of operators $\Sigma$ in the set $\Lambda^{r+1}$ such that $\Sigma(w)=v$.
\label{prop:amongspecial}
\end{proposition}

We defer the verification of both proposition to the following section in favor of proving the main theorem which follows easily from them.  We will write $\Sigma$ for a sequence of operators, $\Sigma^{-1}$ for the same sequence taken in the opposite order, and $\Sigma_2 \Sigma_1$ for the sequence of operators obtained from sequences $\Sigma_1$ and $\Sigma_2$ by applying the operators in $\Sigma_1$ first followed by the ones from $\Sigma_2$.

\begin{proof} Since $\tilde{T}_1=MT(T_1,U)$ for some set of non-core open cycles $U$, the definition of somewhat special implies that $S(T_1)=S(\tilde{T}_1)$.  According to Proposition \ref{prop:makespecial}, there are sequences of operators $\Sigma_1$ and $\Sigma_2$ in $ \Lambda^{r+1}$ such that
    \begin{enumerate}
        \item $\Sigma_1(w)$ has left tableau $S(T_1)$, and
        \item $\Sigma_2(v)$ has left tableau $S(\tilde{T}_1)$
    \end{enumerate}
Since $S(T_1)=S(\tilde{T}_1)$, Proposition \ref{prop:amongspecial} implies that there is another sequence $\Sigma_3$ in  $\Lambda^{r+1}$ such that $ \Sigma_3(\Sigma_1(w))=\Sigma_2(v)$.  In other words,
$$\Sigma^{-1}_2 \Sigma_3 \Sigma_1(w) = v.$$
\end{proof}

Our main approach to the proofs of Propositions \ref{prop:amongspecial} and \ref{prop:makespecial} will be inductive.  The following lemma in a sense serves as the base case.  We show that family of operators $\Lambda^{r+1}$ is enough to generate combinatorial left cells in two special cases. In what follows, we will identify group elements with their corresponding  tableau pairs and apply the operators in $\Lambda^{r+1}$ directly to them via the action described in Section \ref{section:operators}.

\begin{lemma} Let $\tilde{\Lambda}^{r+1}$ consist of the operators $ \Lambda^{r+1} \setminus \{T^{r+1}_{SC}\}.$  Then
\begin{enumerate}
    \item $\tilde{\Lambda}^{r+1}$ generates the irreducible, and consequently the reducible combinatorial left cells of rank $r$ in $W_{n}$ when $n\leq r+1$, and
    \item $\Lambda^{r+1}$ generates the reducible combinatorial left cells of rank $r$ in $W_{r+2}$.
\end{enumerate}

\label{lemma:full}
\end{lemma}

\begin{proof}
The first part is just a restatement of Theorem \ref{theorem:bi}, proved by C.~Bonnaf{\'e} and L.~Iancu.  So consider $w,v$ in $W_{r+2}$ and let $$G_r(w)=(T_1,T_2) \text{ and }G_r(v)=(\tilde{T}_1,\tilde{T}_2)$$ so that $\tilde{T}_1 = MT(T_1,U)$ for some set $U$ of non-core open cycles in $T_1$.  We would like to define a sequence $\Sigma$ of operators in $\Lambda^{r+1}$ so that $\Sigma(w)=v$, or identifying group elements with tableau pairs,
$\Sigma(T_1,T_2)=(\tilde{T}_1,\tilde{T}_2).$

Suppose first that $T_1=\tilde{T}_1$ and the domino with label $r+2$ appears in the same position in $T_2$ and $\tilde{T}_2$. Then by the first part of the lemma, there is a sequence $\Sigma_1 \subset \Lambda^{r+1}$ so that $\Sigma_1(T_1,T_2)'=(\tilde{T}_1,\tilde{T}_2)'.$  But because $T_1=\tilde{T}_1$, this implies $\Sigma_1(T_1,T_2)=(\tilde{T}_1,\tilde{T}_2).$

Next, suppose that $T_1=\tilde{T}_1$ and the domino with label $r+2$ appears in different positions in $T_2$ and $\tilde{T}_2$.  We will say that a tableau of rank $r$ is full if all of its squares of the form $S_{i,r+3-i}$ are occupied.  First suppose that $T_2$ is not full.  We will show that there are sequences $\Sigma'_2$ and $\Sigma''_2$ so that $\Sigma'_2(T_1,T_2)$ and $\Sigma''_2(\tilde{T}_1,\tilde{T}_2)$
have the $r+2$ domino in the same position;  because of what we have already proved, this will be sufficient.  Because $T_2$ and $\tilde{T}_2$ are not full, either their top rows or first columns must have two dominos.  Without loss of generality, assume that it is the top rows, and furthermore, assume that $r+2$ is not in the top row of $T_2,$ for otherwise we could simply swap $(T_1,T_2)$ and $(\tilde{T}_1,\tilde{T}_2)$.  By the first part of the lemma, there is a sequence $\Sigma'_2$ so that $\Sigma'_2(T_1,T_2)'$ has the dominos $r$ and $r+1$ in the top row.   However, note that $T_{\alpha_{r+1} \alpha_{r+2}}$ now can swap the $r+1$ and $r+2$ dominos in the right tableau of $\Sigma'_2(T_1,T_2)$.  If $r+2$ is in the top row of $\tilde{T_2}$, then we are done.  If not, we can repeat the above procedure with $(\tilde{T}_1,\tilde{T}_2)$ to form a sequence $\Sigma''_2.$

Next, suppose that $T_1=\tilde{T}_1$, the domino with label $r+2$ appears in different positions in $T_2$ and $\tilde{T}_2$, and $T_2$ is full.  Let $\Sigma'_3$ be the sequence of operators in $\Lambda^{r+1}$ which arranges the entries of $w$ in decreasing order of absolute values; because $T_2$ is full, the sequences of positive and negative entries in $w$ are both decreasing in absolute value and so such a $\Sigma'_3$ always exists.  By construction,  $\Sigma'_3(T_1,T_2)$ lies in the domain of $T_{SC}$.  If $T^{r+1}_{I\!N}$ is not an element of the sequence $\Sigma'_3$, then the left tableau of $\Sigma'_3(T_1,T_2)$ is $T_1$ and the sought-after sequence of operators is $T_{SC}\Sigma'_3.$  If $T^{r+1}_{I\!N}$ is  in $\Sigma'_3$, then the left tableau of $\Sigma'_3(T_1,T_2)$ differs from  $T_1$ by moving through the cycle $\{r+2\}$.  Then, it is always possible to find a sequence $\Sigma''_3 \subset \Lambda^r$ so that $T_{SC}\Sigma'_3 (T_1,T_2)$ is in the domain of $T^{r+1}_{I\!N}$, perhaps after interchanging the roles of $w$ and $v$.
The desired sequence is then $T^{r+1}_{I\!N} \Sigma''_3 T_{SC} \Sigma'_3.$

Finally, suppose $T_1 \neq \tilde{T}_1$.  The only possible non-core cycle in $T_1$ is $\{r+2\}$ so that $D(r+2,\tilde{T}_1)=D'(r+2, T_1)$.  Because of what we have already proved, it will be enough to find a sequence $\Sigma_4$ so that $\Sigma_4(T_1,T_2)$ and $(\tilde{T}_1,\tilde{T}_2)$ are of the same shape.  Note that because $T_1$ and $\tilde{T}_1$ are of different shapes, either $w$ or $v$ must contain some entries with opposite signs. Without loss of generality, assume that it is $w$;  then it is possible to find a sequence $\Sigma'_4 \subset \Lambda^r$ so that $\Sigma'_4 (T_1,T_2)$ is in the domain of $T^{r+1}_{I\!N}$.  The desired sequence is then $\Sigma_4=T^{r+1}_{I\!N} \Sigma'_4.$

\end{proof}

\section{Somewhat Special Cells}

The goal of this section is to verify Proposition \ref{prop:amongspecial} and prove a few auxiliary facts which we will need for Proposition \ref{prop:makespecial}, whose proof we defer to the next section.
We will say that a sequence $\Sigma$ in  $\Lambda^{r+1}$ is inductive if it does not contain the operators $T_{\alpha_{n-1}\alpha_n}$ and $T_{\alpha_{n}\alpha_{n-1}},$ or the operators $T^{n-1}_{SC}$ and $T^{n-1}_{I\!N}.$
The first lemma describes what happens to the maximal domino in the right tableau of $\Sigma(w)$ for an inductive sequence $\Sigma$ in  $\Lambda^{r+1}$ and is essentially  \cite{garfinkle3}(3.2.7).   While the statement of the result is the same, minor adaptations in the proof are necessary in this more general setting.  We relate the entire proof for completeness.

\begin{lemma} Consider an inductive sequence $\Sigma$ of operators in $ \Lambda^{r+1}$, and let $G_r(w) = (T_1,T_2)$ and $G_r(\Sigma(w))= (\tilde{T}_1,\tilde{T}_2).$
Then:
    \begin{enumerate}
        \item $D(n,\tilde{T}_2)$ is either $D(n,T_2)$ or its image under moving through.
        \item If $\tilde{T}_2(n-1)$ and $T_2(n-1)$ have the same shape, then $D(n,\tilde{T}_2)= D(n,T_2).$
        \item If $D(n,T_2)$ is boxed and $\tilde{T}_2(n-1)$ is somewhat special, then so is $\tilde{T}_2.$
    \end{enumerate}
\label{lemma:technical1}
\end{lemma}

\begin{proof}
The proof of part one proceeds by induction on the size of $\Sigma$.  If $\Sigma$ consists of one operator, then because  $\Sigma$ is inductive the operators $T_{\alpha \beta}$, $T^{r+1}_{SC}$ and $T^{k}_{I\!N}$ for $k\leq r$ (when inductive) leave $D(n,T_2)$ untouched.  The operator $T^{r+1}_{I\!N}$, when inductive, changes $D(n,T_2)$ into $D'(n,T_2)$ iff $n$ is in the extended cycle affected by its action.
The inductive step is an immediate consequence of the fact that $n$ is the maximal entry in $T_2$ and that moving through is an involution.

The other two parts rely on the following consequence of the condition $D(n,\tilde{T}_2) \neq D(n,T_2).$  If this is the case, then either
    \begin{enumerate}
        \item for some $k<n$ the positions of the dominos $D(k,\tilde{T}_2)$ and $D(n,T_2)$ intersect, so that $k$ is in the same cycle as $n$, or
            \label{casei}
        \item for some $k<n$ there are non-core open cycles $c,d$ in $\tilde{T}_2$ and a cycle $e$ in $\tilde{T}_1$ with
$S_b(e)=S_b(d) \in D(k, \tilde{T}_2)$ and
 $S_f(e)=S_b(c) \in D(n, T_2)$,
            implying that $k$ is in the same extended cycle as $n$.
            \label{caseii}
    \end{enumerate}

Both possibilities imply that the shapes of $\tilde{T}_2(n-1)$ and $T_2(n-1)$ are not the same, which proves the second part of the lemma.  For the third part, we will show that $D(n,\tilde{T}_2)=D(n,T_2)$, which will imply that $\tilde{T}_2$ is somewhat special.  We will argue by contradiction, so assume that $D(n,\tilde{T}_2)=D'(n,T_2)$ instead.  If (\ref{casei}) holds, then the boxing condition implies that $D(k,\tilde{T}_2) \cap D(n,T_2)$ is a filled corner in $\tilde{T}_2(n-1)$, contradicting the hypothesis that $\tilde{T}_2(n-1)$ is somewhat special.  On the other hand, if (\ref{caseii}) holds then $D(n,\tilde{T}_2)$ is unboxed as is its entire extended cycle in $\tilde{T}_2$. Then  $S_b(d) \in D(k, \tilde{T}_2)$ is a filled corner  $\tilde{T}_2(n-1)$, contradicting the hypothesis that $T'_2(n-1)$ is somewhat special.
\end{proof}

 The next lemma shows that inductive sequences behave well with respect to domino insertion.
\begin{lemma}
Let $\Sigma$ in $\Lambda^{r+1}$ be an inductive sequence and let $(T_1,T_2)'$ be the tableau pair obtained by deleting the highest-numbered domino from $T_2$ and reversing one step of the insertion procedure for  $T_1$.  Then
$$\Sigma((T_1,T_2)') = (\Sigma(T_1,T_2))'.$$
\label{lemma:induction}
\end{lemma}

\begin{proof}
It is enough to verify this statement when $\Sigma$ consists of just one operator.  When the operator is of the form $T_{\alpha \beta}$, $T^k_{I\!N}$ for $k \leq r$, and $T^{r+1}_{SC}$ for $r+2 < n$ the result is clear.  When $r+2 > n$, $T^{r+1}_{SC}$ is not defined, and when $r+2=n$, it is excluded.  It remains to check the lemma for the operator $T^{r+1}_{I\!N}$ under the assumption that $r+2<n$.

If the action of the operator $T^{r+1}_{I\!N}$ merely swaps two not adjacent dominos in $T_2$, then the result is again clear.  If not, then
 its action on $T_2$ either swaps two domino configurations and then moves through the extended cycle through $r+2$, or vice versa.  Write $T$ for the tableau to which the moving through will be applied.
By Remark \ref{remark:ec}, the extended cycle through $r+2$ in $T$ with respect to $T_1$ contains only non-core cycles.  Thus to show the lemma, it is only necessary to verify the following:

\begin{lemma}  Consider an extended cycle $c$ in $T_2$ relative to $T_1$ which contains only non-core cycles and excludes the cycle $\{n\}$, if one exists.  Write $(\bar{T}_1,\bar{T}_2)$ for the pair  $MT((T_1,T_2), c)$, and form the set $c'$ by deleting the label $n$ from any cycle in $c$.  Write $(\bar{T}'_1,\bar{T}'_2)=(\bar{T}_1,\bar{T}_2)'.$
Then
    \begin{enumerate}
        \item $c'$ is an extended cycle in ${T'_2}$ relative to ${T'_1}$, and
        \item $(\bar{T}'_1,\bar{T}'_2) = MT((T_1,T_2)',c').$
    \end{enumerate}

\end{lemma}
Because $c$ contains only non-core cycles, \cite{garfinkle2}(2.2.9) holds, and the proof of the lemma is identical to that of \cite{garfinkle2}(2.3.3).
\end{proof}

We will need one more lemma to prove Proposition \ref{prop:amongspecial}, which we state presently and whose proof we defer until the next section.  For a somewhat special tableau, it constructs a sequence of operators whose action moves the domino with label $n$ to any other pre-prescribed extremal position.

\begin{lemma}
Consider $(T_1,T_2)$ of rank $r$ with $T_2$ somewhat special.  If $P'$ is a removable pair of adjoining squares in $T_2$, then there is a sequence of operators $\Sigma \subset \Lambda_n$ such that the $n$ domino of the right tableau of $\Sigma(T_1,T_2)$ is in position $P'$ and the left tableau of $\Sigma(T_1,T_2)$ remains equal to $T_1$.
\label{lemma:technical2}
\end{lemma}

Armed with this, we are ready to prove Proposition \ref{prop:amongspecial}, which we restate for the reader's convenience.

\begin{prop:amongspecial}Suppose that $w,v \in W_n$ with $G_r(w)=(T_1, T_2)$ and $G_r(v) = (\tilde{T}_1, \tilde{T}_2)$.  If $T_1 = \tilde{T}_1$ and both are somewhat special, then there is a sequence of operators $\Sigma$ in the set $\Lambda^{r+1}$ such that $\Sigma(w)=v$.
\end{prop:amongspecial}

\begin{proof} We argue by induction and assume that Theorem \ref{theorem:generators} holds for numbers smaller than $n$.  Let $P'$ be the position of the $n$ domino in $\tilde{T}_2$.  By Lemma \ref{lemma:technical2}, there is a sequence $\Sigma_1$ of operators so that the right tableau of $\Sigma_1(T_1,T_2)=(T^1_1,T^1_2)$ has its $n$ domino in position $P'$ and $T_1=T^1_1$.  Note that the left tableaux of $(T^1_1,T^1_2)'$ and $(\tilde{T}_1, \tilde{T}_2)'$ must be the same since $\tilde{T}_2$ and $T^1_2$ have their $n$ dominos in the same position, and by hypothesis, $T^1_1=\tilde{T}_1$.  By induction, we can construct an inductive sequence $\Sigma_2$ such that $\Sigma_2 (T^1_1,T^1_2)' = (\tilde{T}_1, \tilde{T}_2)'$, which, by Lemma \ref{lemma:induction}, also equals $(\Sigma_2 (T^1_1,T^1_2))'$.

We would like to show that $\Sigma_2 (T^1_1,T^1_2)= (\tilde{T}_1, \tilde{T}_2)$, which will complete the proof.  Since $(T^1_1,T^1_2)'$ and $(\tilde{T}_1, \tilde{T}_2)'$
must have the same shape, Lemma \ref{lemma:technical1} implies that the $n$ dominos of the right tableau of $\Sigma_2 (T^1_1,T^1_2)$ and $T^1_2$ must be in the same position, which also happens to be the position of the $n$ domino in  $\tilde{T}_2$.  But since $\Sigma_2 (T^1_1,T^1_2)'= (\tilde{T}_1, \tilde{T}_2)'$, the right tableaux of $\Sigma_2 (T^1_1,T^1_2)$ and  $(\tilde{T}_1, \tilde{T}_2)$ must agree.  Hence their left tableaux must be of the same shape, and because $T_1^1=\tilde{T_1}$, the rest follows.

\end{proof}

\section{Technical Lemmas}

In this section, we verify Proposition \ref{prop:makespecial} and Lemma \ref{lemma:technical2}.  They are both adaptations of lemmas from \cite{garfinkle2} and follow from case by case analyses of relative domino positions.

\begin{prop:makespecial} Given $w\in W_n$, there is a sequence of operators $\Sigma$ in  $\Lambda^{r+1}$ such that the left tableau of $\Sigma(w)$ is somewhat special.
\end{prop:makespecial}

\begin{proof}
Consider $w \in W_n$ and let $G_r(w) = (T_1,T_2)$.  Our goal is to find a sequence $\Sigma \in \Lambda^{r+1}$ so that $\Sigma(w)$ has a somewhat special left tableau.
We will superficially follow the outline of the proof of \cite{garfinkle3}(3.2.4), but the new situations in the general rank case will require a slightly different approach and our cases are somewhat different. We will use induction on $n$, assuming both, that Proposition \ref{prop:makespecial} and Theorem \ref{theorem:generators} are true for values smaller than $n$.

Let $(\bar{T}_1,\bar{T}_2) = (T_1,T_2)'$.  By induction, we know that there is an inductive sequence $\Sigma_1$  so that the left tableau of $ \Sigma_1((T_1,T_2)')= (\bar{T}^1_1,\bar{T}^1_2)$
is somewhat special.  Let $\Sigma_1(T_1,T_2)=({T}^1_1,{T}^1_2)$.  Then Lemma \ref{lemma:induction} implies $({T}^1_1,{T}^1_2)'= (\bar{T}^1_1,\bar{T}^1_2)$.  There are two cases:
    \begin{enumerate}
        \item If ${T}^1_2$ is somewhat special, then so is $T^1_1$ and we are done.
        \item If ${T}^1_2$ is not somewhat special, then one of its non-core open cycles is not boxed.  But since $\bar{T}^1_2$ itself is somewhat special, the only possibility is that $\{n\}$ is an unboxed extended open cycle in $T^1_2$ relative to $T^1_1$.
    \end{enumerate}
We assume the latter is true and without loss of generality, take the domino $D(n,T^1_2)$ to be horizontal writing $\{S_{ij},S_{ij+1}\}$ for its underlying squares.  Because $\{n\}$ is a cycle, $D(n,T^1_2)$ is unboxed, and $T^1_2$ is not somewhat special, $S_{ij+1}$ must be a filled corner and $S_{i+1,j}$ must be empty hole in $T_1^1$ implying that $S_{i,j+1}$ is a square of type $X$.
We will examine a number of cases in order to find a sequence of operators which will take $T^1_2$ to a somewhat special tableau.

Our general goal will be to find two pairs of adjacent squares $P_1$ and $P_2$ in $T^1_2$, so that $P_1$ is extremal in $T^1_2(n-1)$ and $P_2$ is extremal in $T^1_2(n-1) \setminus P_1.$  If this is possible, then induction on Theorem \ref{theorem:generators} provides a sequence $\Sigma_2$ so that
$(\bar{T}^2_1,\bar{T}^2_2)=\Sigma_2(\bar{T}^1_1,\bar{T}^1_2)$ with $\bar{T}^2_1=\bar{T}^1_1$ and the $n-1$ domino of $\bar{T}^2_2$ is $P_1$ while its $n-2$ domino is $P_2$.  Let $({T}^2_1,{T}^2_2) = \Sigma_2({T}^1_1,{T}^1_2)$.  If $({T}^2_1,{T}^2_2)$ is in the domain of either $T_{\alpha \beta}= T_{\alpha_{n-1} \alpha_n}$ or $T_{\alpha_{n} \alpha_{n-1}}$, then let $({T}^3_1,{T}^3_2)=T_{\alpha \beta}({T}^2_1,{T}^2_2).$  Again by induction, it is then possible to find a sequence $\Sigma_4$ so that
$(\bar{T}^4_1,\bar{T}^4_2) = \Sigma_4({T}^3_1,{T}^3_2)'$ is somewhat special.  Finally, let $$({T}^4_1,{T}^4_2) = \Sigma_4({T}^3_1,{T}^3_2).$$
This turns out to be somewhat special, and $\Sigma = \Sigma_4 (T_{\alpha \beta}) \Sigma_2 \Sigma_1$ is the promised sequence of operators.

If pairs of adjacent squares $P_1$ and $P_2$ cannot be suitably chosen, we will find that $T^1_2$ satisfies the hypotheses of Lemma \ref{lemma:full}, which provides the desired sequence of operators.  To simplify the details, we will let a number of our cases overlap.

{\it Case 1.}\hspace{.05in}
Assume $i>1$, which, because of our choice of boxing, implies that $i \geq 3$. Let $s$ and $t$ be the lengths of the $i-1$st and $i-2$nd rows of $T^1_2$, and further assume that $t-s$ and $s-j$ are both greater than $1$.  Define $P_1=\{S_{i-2,t-1},S_{i-2,t}\}$ and $P_2=\{S_{i-1,s-1},S_{i-1,s}\}$.   Then $({T}^2_1,{T}^2_2)$ lies in the domain of $T_{\alpha_{n-1} \alpha_n}$, which acts on $T^2_2$ by interchanging the $n$ and $n-1$ dominos.  That  $({T}^4_1,{T}^4_2)$ is somewhat special follows from Lemma \ref{lemma:technical1} since $\bar{T}^4_2$ is somewhat special and  $P_1$ is boxed in $\bar{T}^1_2$ which is also somewhat special.

{\it Case 2.}\hspace{.05in}  Assume the same situation as in Case 1, but suppose $t=s > j+1$.  The proof is identical to Case 1 if we let $P_1=\{S_{i-2,t},S_{i-1,t}\}$ and $P_2=\{S_{i-2,t-1},S_{i-1,t-1}\}$.

{\it Case 3.}\hspace{.05in}  Assume the same situation as in Case 1, but suppose $t=s = j+1$.  Define $P_1=\{S_{i-2,t},S_{i-1,t}\}$ and $P_2=\{S_{i-2,t-1},S_{i-1,t-1}\}$.  The action of $T_{\alpha_{n-1} \alpha_n}$ on $T^2_2$ this time is not just an interchange and is described by case 1 of \ref{section:alphabeta}.  That $T^4_1$ is somewhat special follows because $\bar{T}^4_2$ is somewhat special and $D(n,T^4_2)$ occupies the squares $\{S_{i-2,j+1},S_{i-1,j+1}\}$ and therefore is boxed.

{\it Case 4.}\hspace{.05in} Assume that $S_{i+2,j-1}$ is  in $T^1_2$.  Because $\bar{T}^1_2$ is somewhat special, this implies $S_{i+3,j-1}$ is also in $T^1_2$, otherwise $\bar{T}^1_2$ would have a filled corner.  Let $l$ be the length the of $j-1$st column of $T^1_2$ and define $P_1=\{S_{l-1,j-1},S_{l,j-1}\}$ $P_2=\{S_{l-3,j-1},S_{l-2,j-1}\}$.  The rest of the proof follows as before, except the operator $T_{\alpha_{n} \alpha_{n-1}}$ must be used.

{\it Case 5.}\hspace{.05in} Assume $S_{i+2,j-1}$ is  not in $T^1_2$, but that $S_{i+2,j-2}$ is.  Let $P_2$ be the pair of squares $\{S_{i,j-1},S_{i+1,j-1}\}$.  The definition of $P_1$ is more involved.  Let $m$ be the smallest positive number so that either $\{S_{i+m+1,j-m-1}\}$ is not a square in $T^1_2$ or $\{S_{i+m+1,j-m}\}$ is a square in $T^1_2$. We examine the former as case (a), the latter as case (b), and if neither occurs as case (c).
    \begin{enumerate}
        \item[(a)] Define $P_1=\{S_{i+m,j-m-1},S_{i+m,j-m}\},$
        \item[(b)] Let $p$ be the length of the $(j-m)$th column of $T^1_2$ and define $P_1$ to be the pair of squares $\{S_{p-1,j-m-1},S_{p,j-m}\}.$
    \end{enumerate}
In either of these two cases, the proof follows as before using the operator $T_{\alpha_{n} \alpha_{n-1}}$.
    \begin{enumerate}
        \item[(c)] In this case, no adequate extremal position for $P_1$ exists below $D(n,T^1_2)$.  If $i=1$, then $T^1_2$ in fact has to satisfy the hypotheses of Lemma \ref{lemma:full}, so we can assume (as described in Case (1)) that $i>2$.  We may further assume that $T^1_2$ does not satisfy the hypotheses of any of the previous cases, so that we can also assume $s=j+1$ and $t=j+2$. Define $m$ to be the smallest positive number so that either $\{S_{i-m-1,j+m+1}\}$ is not a square in $T^1_2$ or $\{S_{i-m,j+m+1}\}$ is a square in $T^1_2$.  As before, define $P_1$ to be respectively $\{S_{i-m-1,j+m},S_{i-m,j+m}\}$ or $\{S_{i-m,q-1},S_{i-m,q}\}$ where  $q$ is the length of the $i-m$th row of $T^1_2$. In either of these two cases, the proof follows as before using the operator $T_{\alpha_{n-1} \alpha_{n}}$. If no such $m$ exists however, then $T^1_2$ again has to satisfy the hypotheses of Lemma \ref{lemma:full}, which provides the necessary sequence of operators.
    \end{enumerate}

{\it Case 6.}\hspace{.05in}  Assume $j>2$ and $S_{i+2,j-2} \notin T^1_2.$  Note that since $S_{i+2,j-2} \notin T^1_2,$  the square $S_{i,j-1}$ cannot lie in the core of $T^1_2.$  Thus let $P_1=\{S_{i+1,j-2},S_{i+1,j-1}\}$ and $P_2=\{S_{i,j-2},S_{i,j-1}\}.$
The proof proceeds as before, with $D(n,T^4_2)=D'(n,T^3_2) = \{S_{i+1,j-1},S_{i+1,j}\}.$

{\it Case 7.}\hspace{.05in} Assume $j=2$ and $S_{i+2,j-1} \notin T^1_2.$  If we assume that $T^1_2$ satisfies the hypotheses of none of the previous cases, then let $P_2=\{S_{i,j-1},S_{i+1,j-1}\}$ and the proof proceeds exactly as in Case 5(c).

{\it Case 8.}\hspace{.05in} Assume $j=1$, so that $i$ cannot equal 1, implying $i>2$.  If we assume that $T^1_2$ satisfies the hypotheses of none of the other cases above, then we can take $s=2$ and $t=3$.  Let $P_2=\{S_{i-1,1},S_{i-1,2}\}.$ Again, the proof proceeds exactly as in Case 5(c).
It is a quick check to see that our cases cover all possibilities, and so the proof of the lemma is complete.
\end{proof}

The following lemma is a general rank version of \cite{garfinkle3}(3.2.9).  We adapt its proof to the general rank case.  There are two complications:  the existence of a nontrivial core and the fact that the tableau is only somewhat special, so we cannot count on all holes being filled.   Almost entirely, cases which do not match those in the rank zero proof will degenerate into the setting of Lemma \ref{lemma:full}, where we proved the main result for tableaux whose rank is large relative to $n$.

\begin{lemma:technical2}
Consider $(T_1,T_2)$ of rank $r$ with $T_2$ somewhat special.  If $P'$ is a removable pair of adjoining squares in $T_2$, then there is a sequence of operators $\Sigma \subset \Lambda_n$ such that the $n$ domino of the right tableau of $\Sigma(T_1,T_2)$ is in position $P'$ and the left tableau of $\Sigma(T_1,T_2)$ remains equal to $T_1$.
\end{lemma:technical2}

\begin{proof}
The proof is by induction where we assume both, that Lemma \ref{lemma:technical2} and Theorem \ref{theorem:generators} are true for values smaller than $n$.  Write $P$ for the set of squares occupied by the $n$ domino in $T_2$ and assume that $P \neq P'$.
Without loss of generality assume that $P$ is horizontal, writing $P=\{S_{ij},S_{i,j+1}\}$.

\vspace{.1in}
{\it Cases 1-4.}\hspace{.05in}
For the first cases, also assume that $P'= \{S_{kl},S_{k,l+1}\}$ is horizontal.  We first define an auxiliary domino position $P_1$.
\begin{enumerate}
    \item
If $k=i-1$, we must have $l \geq j+2$ since $P'$ is extremal. If $S_{i-1,j}$ is in the core of $T_2$, then $l \geq j+3.$
In either case, let $P_1=\{S_{i-1,l-2},S_{i-1,l-1}\}.$

\item  If $k<i-1$, let $u$ be the length of the $(i-1)$st row of $T_2$.  Then $l \geq u+1$ and we can set $P_1=\{S_{i-1,u-1},S_{i-1,u}\}.$

\item If $k=i+1$, then $l\leq j-2$.  Let $P_1=\{S_{i,j-2},S_{i,j-1}\}.$  Note that if $S_{i,j-2}$ is in the core of $T_2$, then $P'$ could not have been extremal.

\item If $k>i+1$, then let $u$ be the length of the $(k-1)$st row of $T_2$.  Let $P_1= \{S_{k-1,u-1},S_{k-1,u}\}$.

\end{enumerate}
Proceeding by induction as in the proof of  Proposition \ref{prop:makespecial}, we can find an inductive sequence $\Sigma_1$  so that $\Sigma_1(T_1,T_2)'$ has the same left tableau as $(T_1,T_2)'$ and whose right tableau has domino $n-1$ in position $P'$, and $n-2$ in position $P_1$.  By Lemma \ref{lemma:technical1}, $\Sigma_1(T_1,T_2)$ has left tableau $T_1$.  The desired sequence is then  $T_{\alpha_{n-1},\alpha_n} \Sigma_1$ in the first two cases and $T_{\alpha_{n},\alpha_{n-1}} \Sigma_1$ in the second two.

\vspace{.1in}
{\it Cases 5-7.}\hspace{.05in}    For the rest of the cases assume that $P'= \{S_{kl},S_{k+1,l}\}$ is vertical.  We again define an auxiliary domino position $P_1$.

\begin{enumerate}
    \item[(5)] If $k+1=i-1$, let $P_1= \{S_{i-1,l-2},S_{i-1,l-1}\}$ as in Case (1).

    \item[(6)]  If $k+1<i-1$, let $P_1=\{S_{i-1,u-1},S_{i-1,u}\} $ as in Case (2).

    \item[(7)] If $k=i+1$ and $l <j-1$, let $P_1=\{S_{i,j-2},S_{i,j-1}\}$ as in Case (3).
\end{enumerate}
The proof proceeds exactly as in the first four cases, with $T_{\alpha_{n-1},\alpha_n}$ applied in Cases 5 and 6, and $T_{\alpha_{n},\alpha_{n-1}}$ in Case 7.

\vspace{.1in}
{\it Case 8.}\hspace{.05in} Again assume that $P'= \{S_{kl},S_{k+1,l}\}$ and suppose that $k>i+1$.
    \begin{enumerate}
        \item[{\it i.}] Suppose that $P''$ is an extremal domino shape in $T_2 \setminus P \setminus P'$.  If $P''$ is above $P'$ but not above $P$, let $P_1=P'$ and proceed as in Case (4).  If $P''$ is below $P'$ and horizontal, then using cases (1)-(4) we can find a sequence so that the resulting right tableau will have the $n$ domino in position $P''$. Then we are in the setting of cases (5)-(7).  If $P''$ is  is below $P'$ and vertical, by induction it is possible to find a sequence so that the resulting right tableau will have the $n-1$ domino in position $P'$ and the $n-2$ domino in position $P''$.  After applying $T_{\alpha_n,\alpha_{n-1}}$, we are in the setting of the transpose of cases (1)-(4).  If $P''$ is above $P$ and also extremal in $T_2$ itself, then arguments similar to the above apply.

        \item[{\it ii.}] Suppose that $P''= \{S_{k+1,l-2},S_{k+1,l-1}\}$ is extremal in $T_2 \setminus P'$ and $P'''= \{S_{k,l-2},S_{k,l-1}\}$ is extremal in  $T_2 \setminus P' \setminus P''.$ By induction it is possible to find a sequence so that the resulting right tableau will have the $n-1$ domino in position $P'$, $n-2$ domino in position $P''$, and the $n-3$ domino in position $P'''.$  After applying $T_{\alpha_{n-2},\alpha_{n-1}}$, we are in the setting of cases (1)-(4).  A similar argument works if
                \begin{enumerate}
                    \item $P''= \{S_{k,l-1},S_{k+1,l-1}\}$ is extremal in $T_2 \setminus P'$ and the domino shape $P'''= \{S_{k-1,l-1},S_{k-1,l}\}$ is extremal in  $T_2 \setminus P' \setminus P''.$
                    \item $P''= \{S_{i-2,j+1},S_{i-1,j+1}\}$ is extremal in $T_2 \setminus P$ and the domino shape $P'''= \{S_{i-2,j},S_{i-1,j}\}$ is extremal in  $T_2 \setminus P \setminus P''.$
                    \item $P''= \{S_{i-1,j},S_{i-1,j+1}\}$ is extremal in $T_2 \setminus P$ and the domino shape $P'''= \{S_{i-1,j-1},S_{i,j-1}\}$ is extremal in  $T_2 \setminus P \setminus P''.$
                \end{enumerate}
    \end{enumerate}

If there is no extremal $P''$ in  $T_2 \setminus P \setminus P'$ above $P'$ but not above $P$, $S_{i+1,j-2} \in T_2$, but $S_{i+1,j-1} \notin T_2$.  It then follows by an easy but tedious inspection that if we are not in any of the above cases, then we must be within the scope of Lemma \ref{lemma:full}.

\vspace{.1in}
{\it Case 9.}\hspace{.05in}  Assume that $P$ and $P'$ are both boxed and that $k=i-1$, so that $P$ and $P'$ intersect in the square $S_{i,j+1}$.  Since they are both boxed, we know that $S_{ij+1}$ must be of type $W$, and consequently, that $i$ is even.
    \begin{enumerate}
        \item[{\it i.}]  Assume that $i>2$, and hence that $i\geq 4$.  Let $u$ be the length of the  row $i-2$ of $T_2$ and let $s$ be the length of row  $i-3$.
                \begin{enumerate}
                    \item If $u=j+1$, let $P_1=\{S_{i-1,j+1},S_{i-2,j+1}\}$ and $P_2=\{S_{i-1,j},S_{i-2,j}\}.$  Note that $S_{i-2,j}$ cannot be in the core of $T_2$, since $D(n,T_2)$ occupies the square $S_{ij}$. By induction, we can find a sequence $\Sigma_1$ so that the $n-1$ domino in the right tableau of $\Sigma_1(T_1,T_2)'$ is $P_1$ and the $n-2$ domino is $P_2$.  Then the $n$ domino of $T_{\alpha_{n-1},\alpha_n}\Sigma_1(T_1,T_2)$ occupies the squares $\{S_{i,j+1},S_{i-1,j+1}\}$, as desired.

                        \item Assume that $s=u>j+1$.  If we let $P''=\{S_{i-3,u},S_{i-2,u}\}$, then Case 6 provides a sequence $\Sigma_1$ so that the right tableau of $\Sigma_1(T_1,T_2)$ has the $n$ domino in position $P''$.  The transpose of the first two cases now provides a sequence $\Sigma_2$ such that the $n$ domino of $\Sigma_2\Sigma_1(T_1,T_2)$ is $P'$, as desired.

                        \item Assume that $s>u>j+1$.  If there is an extremal domino shape $P''$ in $T_2$ strictly above $P$, we can use the already considered cases (or their transposes), to find a sequence $\Sigma_1$ so that $\Sigma_1(T_1,T_2)$ has left tableau $T_1$ and the $n$ domino is in position $P''$ in its right tableau. Then, again by the already considered cases (or their transposes), we can find a sequence $\Sigma_2$ so that $\Sigma_2\Sigma_1(T_1,T_2)$ has left tableau $T_1$ and the $n$ domino is in position $P'$ in its right tableau, as desired. A similar argument works if there is an extremal domino shape $P''$ in $T_2$ strictly below $P$.  If neither is the case, then we are in the scope of Lemma \ref{lemma:full}, which provides the necessary sequence.
                    \end{enumerate}
         \item[{\it ii.}] If $i=2$ and $j>2$, then this case is parallel to the previous one.  If $i=2$ and $j=2$, then we are either in the setting of Lemma \ref{lemma:full} or Case 8.  Finally, if $i=2$ and $j=1$, we are in rank zero, which has been considered in \cite{garfinkle3}.
    \end{enumerate}

\vspace{.1in}
{\it Case 10.}\hspace{.05in}   Assume that $P$ and $P'$ are both boxed and that $k=i+1$ and $l=j-1$.  The boxing condition implies that $S_{i,j+1}$ must be of type $W$, and consequently, that $i$ is even.  Let $r$ be the length of the $i-1$st row of $T_2$.  If $u=j+1$, then let $P''=\{S_{i-1,j+1},S_{i,j+1}\}$.  By Case 9, we can find a sequence $\Sigma_1$ so that $\Sigma_1(T_1,T_2)$ has left tableau $T_1$ and its right tableau has the $n$ domino in position $P''$.  Then Case 2 provides a sequence $\Sigma_2$ so that $\Sigma_2\Sigma_1 (T_1,T_2)$ has left tableau $T_1$ and its right tableau has the $n$ domino in position $P'$. If $u>j+1$, then because $T_2$ is somewhat special it cannot have filled corners, so that $u \geq j+3$.  But then we can let $P''=\{S_{i-1,r-1},S_{i,r}\}$ and the proof is the same as when $r=j+1$.

\vspace{.1in}
{\it Case 11.}\hspace{.05in}   Assume that $P'$ is boxed and $P$ is not.  Our goal is to find a sequence of operators which will move the $n$ domino to a boxed position, reducing our work to the previous cases.  Because $T_2$ is somewhat special, it cannot have filled corners, so that $S_{ij}$ must be of type $W$ and consequently, $i\geq 2$.  Let $u$ be the length of the $i-1$ row of $T_2$.  If $u=j+1$, let $\Sigma_1$ be the sequence of operators such that the right tableau of $\Sigma_1 (T_1,T_2)'$ is somewhat special; it exists by induction and Proposition \ref{prop:makespecial}.  The shape of the tableaux $\Sigma_1(T_1,T_2)'$ is $(shape(T_2) \setminus S_{i-1,j+1}) \cup S_{ij}$.  If $T^1_2$ is the right tableau of $\Sigma_1(T_1,T_2)$, then Lemma \ref{lemma:technical1}, implies that $D(n,T^1_2)=D'(n,T_1) = \{S_{i-1,j+1},S_{i,j+1}\}$, which is boxed, as desired.
If $u>j+1$, then because $T_2$ is somewhat special and cannot have filled corners, $u \geq j+3$.  Let $P''=\{S_{i-1,u-1},S_{i-1,u}\}$.  Then from Cases 1 through 4, where a boxing condition was not assumed, we can find a sequence of operators $\Sigma_1$ so that the $n$ domino in the right tableau of $\Sigma_1(T_1,T_2)$ is $P''$, which, because $T_2$ is special, must be boxed.

\vspace{.1in}
{\it Case 12.}\hspace{.05in} Assume that $P'$ is unboxed in $T_2$.  Since Cases 1 through 8 were verified without any boxing assumptions, we only have to check this case when
    \begin{enumerate}
        \item[{\it i.}] $k=i-1$ so that $P$ and $P'$ intersect in the square $S_{i,j+1}$, or
        \item[{\it ii.}] $k=i+1$ and $l=j-1$.
    \end{enumerate}
Furthermore, we can also take $P$ to be boxed, as otherwise the boxing conditions imply in both cases that $S_{i,j+1}$ is a filled corner, contradicting the assumption that $T_2$ is somewhat special.  With this assumption, the transpose of (i) is the last case considered in  \cite{garfinkle3}(3.2.9), and we omit its proof.  So assume that we are in case (ii).  Then the type of $S_{ij}$ is $X$.
    \begin{enumerate}
        \item[(a)] Assume that there are no extremal domino positions strictly below $P'$ or above $P$.  If $i>1$ and $S_{i-1,j+2}$ is not a square in $T_2$, or if $S_{i+2,j-2} \in T_2$ and $S_{i+3,j-2}\notin T_2$, then we can use case (i) above to find a sequence $\Sigma_1$ so that the right tableau of $\Sigma_1(T_1,T_2)$ has the $n$ domino in position $\{S_{i,j+1},S_{i-1,j+1}\}$.  Then $\Sigma_1(T_1,T_2)$ is in the setting of Case 2.  If neither of the two possibilities above is true, then we are in the setting of Lemma \ref{lemma:full}.
        \item[(b)] If there is an extremal domino position $P''$ strictly below $P'$, then using induction we can find a sequence $\Sigma_1$ so that the right tableau of $\Sigma_1(T_1,T_2)$ has domino $n-1$ in position $P''$ and $n-2$ in position $P'$.  Then $T_{\alpha_n \alpha_{n-1}}\Sigma_1 (T_1,T_2)$ is within scope of the previous cases.  If there is an extremal domino position $P'$  strictly  above $P$, then by the previous cases, there is a sequence $\Sigma_1$ so that the $n$ domino of  the right tableau of  $\Sigma_1(T_1,T_2)$ is $P''$.  We can use induction to find a sequence $\Sigma_2$ which puts the $n-2$ domino in position $P$ and the $n-1$ domino in position $P'$.  Then the desired sequence is $T_{\alpha_n \alpha_{n-1}}\Sigma_2 \Sigma_1 (T_1,T_2).$
    \end{enumerate}

\end{proof}


\begin{thebibliography}{10}


\bibitem{bonnafe:iancu}
C.~Bonnaf{\'e} and L.~Iancu.
\newblock Left cells in type $B\sb n$ with unequal parameters.
\newblock {\em Represent. Theory}, 7:587--609, 2003.

\bibitem{bgil}
C.~Bonnaf{\'e}, M.~Geck, L.~Iancu, and T.~Lam.
\newblock On domino insertion and {K}azhdan--{L}usztig cells in type
$B_n$, {\em Progress in Math. (Lusztig Birthday Volume).}
\newblock Birkhauser, to appear.


\bibitem{garfinkle1}
D.~Garfinkle.
\newblock On the classification of primitive ideals for complex classical {L}ie
  algebras (I).
\newblock {\em Compositio Math.}, 75(2):135--169, 1990.

\bibitem{garfinkle2}
D.~Garfinkle.
\newblock On the classification of primitive ideals for complex classical {L}ie
  algebras (II).
\newblock {\em Compositio Math.}, 81(3):307--336, 1992.

\bibitem{garfinkle3}
D.~Garfinkle.
\newblock On the classification of primitive ideals for complex classical {L}ie
  algebras (III).
\newblock {\em Compositio Math.}, 88:187--234, 1993.

\bibitem{geck:leftcells}
M.~Geck.
\newblock Constructible characters, leading coeffcients and left cells for finite {C}oxeter groups
with unequal parameters.
\newblock {\em Represent. Theory}, 6:1-30, 2002.


\bibitem{geck:pfeiffer}
M.~Geck and G.~Pfeiffer.
\newblock Characters of finite {C}oxeter groups and {I}wahori {H}ecke algebras,
\newblock {\em London Math. Soc. Monographs}, New Series 21, Oxford University Press, 2000.

\bibitem{gordon:quiver}
I.~G.~Gordon  \newblock Quiver varieties, category $\mathcal{O}$
 for rational {C}herednik algebras, and {H}ecke algebras,
{\tt arXiv:math.RT/0703150}.

\bibitem{gordon:calogero}
I.~G.~Gordon and M.~Martino. \newblock {C}alogero-{M}oser Space, reduced
rational {C}herednik algeras and two-sided cells,
{\tt arXiv:math.RT/0703153}.

\bibitem{hopkins:domino}
B.~Hopkins.
\newblock Domino tableaux and single-valued wall-crossing operators.  Ph.D. dissretation, University of Washington, 1997.

\bibitem{knuth:art3}
{D.~E.~Knuth,} {The art of computer programming. {V}olume 3},
Addison-Wesley Publishing Co., Reading, Mass.-London-Don
              Mills, Ont., 1973.


\bibitem{lusztig:unequal}
G.~Lusztig.
\newblock {\em {H}ecke algebras with unequal parameters}, volume~18 of {\em CRM
  Monograph Series}.
\newblock American Mathematical Society.


\bibitem{mcgovern:ssmap}
W.~M.~McGovern. \newblock On the {S}paltenstein-{S}teinberg map for
classical {L}ie algebras.
\newblock{\em Comm. Algebra}, 27(6):2979--2993, 1999.



\bibitem{pietraho:components}
T.~Pietraho. \newblock Components of the {S}pringer fiber and domino
tableaux. \newblock{\em Journal of Algebra}, 272 (2):711--729, 2004.

\bibitem{pietraho:rscore}
T.~Pietraho.
\newblock A relation for domino {R}obinson-{S}chensted algorithms,
{\em Annals of Combinatorics}, to appear.

\bibitem{pietraho:equivalence}
T.~Pietraho.
\newblock Equivalence classes in the {W}eyl groups of type $B_n$.,
{\em Journal of Algebraic Combinatorics}, 27(2):247-262, 2008.

\bibitem{pietraho:constructible}
T.~Pietraho.
\newblock Cells and constructible representations in type $B_n$.,
{\tt math.CO/0710.3846}.

\bibitem{stanton:white}
Dennis~W. Stanton and Dennis~E. White.
\newblock A {S}chensted algorithm for rim hook tableaux.
\newblock {\em J. Combin. Theory Ser. A}, 40(2):211--247, 1985.

\bibitem{taskin:plactic}
M.~Taskin.
\newblock Plactic relations for $r$-domino tableaux.
{\tt math.RT:0803.1148}.

\bibitem{vanleeuwen:rank}
M.~A.~A. van Leeuwen.
\newblock The {R}obinson-{S}chensted and {S}chutzenberger algorithms, an elementary
  approach.
\newblock {\em Electronic Journal of Combinatorics}, 3(2), 1996.

\bibitem{vogan:tau}
D. Vogan.
\newblock A generalized $\tau $-invariant for the primitive spectrum of a
   semisimple {L}ie algebra
\newblock {\em Math. Ann.} 242(3):209-224, 1979.


\end{thebibliography}
\end{document}